\documentclass[11pt]{article}
\usepackage{amsmath, amssymb, amsfonts}
\usepackage{amsthm}
\usepackage{geometry}
\usepackage{graphicx}
\usepackage{bm}
\usepackage{booktabs}
\usepackage{cleveref}
\usepackage{tikz}
\usetikzlibrary{arrows.meta, decorations.pathreplacing}

\newcommand{\norm}[1]{\left\|#1\right\|}

\newcommand{\dt}{\Delta t}
\newcommand{\R}{\mathbb{R}}
\newcommand{\bdf}[1]{\delta_q #1}
\newcommand{\ds}{\partial_s}
\newcommand{\abs}[1]{\left|#1\right|}

\newtheorem{theorem}{Theorem}
\newtheorem{corollary}{Corollary}
\newtheorem{assumption}{Assumption}
\newtheorem{remark}{Remark}
\newtheorem{lemma}{Lemma}

\title{High-Order Unconditionally Stable IMEX BDF--FEM Schemes for the Hodgkin--Huxley Model on a Branched Axon}
\author{Hayden Reed and Qingguang Guan\thanks{\raggedright Corresponding author: Qingguang Guan.\\
Email addresses: \texttt{Hayden.Reed@usm.edu} (Hayden Reed),
\texttt{Qingguang.Guan@usm.edu} (Qingguang Guan).}\\[0.5em]
\small School of Mathematics and Natural Sciences, The University of Southern Mississippi\\
\small Hattiesburg, MS 39406, USA}
\date{\today}

\usepackage{etoolbox}
\usepackage{placeins}

\pretocmd{\section}{\FloatBarrier}{}{}
\patchcmd{\thebibliography}{\sloppy}{\setlength{\itemsep}{0pt}\setlength{\parsep}{0pt}\sloppy}{}{}

\begin{document}

	\maketitle

\begin{abstract}
We construct and analyze high-order implicit--explicit backward
differentiation formula finite element (IMEX BDF--FEM) schemes for the
Hodgkin--Huxley cable system on a branched metric graph.  The linearly
implicit schemes combine constant-step BDF methods of orders
$q=1,\ldots,5$ with conforming finite elements of any fixed degree $p\ge1$.
Extrapolated gate variables are smoothly masked before evaluating the
conductances, preserving coercivity without projecting the numerical
gates.  A voltage cutoff makes the gate rates globally bounded and
Lipschitz.  Each timestep requires one linear voltage solve followed by
independent nodal gate updates.  Under exact spatial integration, we prove
unconditional stability of the complete fully discrete regularized scheme
in discrete history energies.  Under regularity and starting-history
assumptions, the voltage and reconstructed gates converge in
$L^2(\Gamma)$ at the optimal rate $O(h^{p+1}+\Delta t^q)$, without a
timestep restriction or an $h$--$\Delta t$ coupling.  The same rate holds
with numerical quadrature when its full residual defect is $O(h^{p+1})$.
Manufactured solutions recover the predicted orders through BDF5--$P_4$.
Action-potential simulations illustrate pulse propagation and splitting,
while a controlled work--precision comparison within the proposed family
shows that higher-order variants can reach a prescribed accuracy at much lower
runtime along the tested joint-refinement paths.
\end{abstract}



	\section{Introduction}

	The Hodgkin--Huxley cable model describes action-potential propagation by
	coupling a parabolic voltage equation to three ionic gate equations
	\cite{hodgkin1952quantitative}.  A branched axon is naturally represented
	as a metric graph: each edge is a cable segment, and voltage continuity
	and Kirchhoff current balance couple the branches at their junctions.
	Rapid voltage fronts and distinct gate time scales make the system stiff.
	Accurate simulation therefore requires spatial resolution together with
	stable time integration and error control for the coupled system.

	Neuronal methods include branched-cable solvers~\cite{hines1984efficient},
	exponential gate updates~\cite{rush1978practical}, and first- and
	second-order exponential integrators that prevent voltage overshoot for
	space-clamped Hodgkin--Huxley-like ODEs~\cite{borgers2013exponential}.
	For the cable PDE--ODE system, Mascagni~\cite{mascagni1990backward}
	proved an $O(\dt+h^2)$ error bound using backward Euler in time and
	centered finite differences in space.  Other cable discretizations address
	high-order spatial accuracy
	\cite{omurtag2010spectral}, positivity and boundedness
	\cite{aderogba2018positive}, and charge conservation
	\cite{ashraf2026finitevolume}.  Branchwise predictor--corrector methods
	decouple cable segments at junctions~\cite{rempe2006predictor}; a locally
	adaptive BDF2 finite-difference method exhibits second-order convergence in
	numerical tests~\cite{kublik2016adaptive}.
	Amsallem and Nordstr\"om~\cite{amsallem2013highorder} established
	well-posedness for the continuous problem and semi-discrete energy
	estimates for high-order SBP--SAT schemes on branched
	Hodgkin--Huxley cables.

	Related neuronal calcium models also require a numerical treatment of
	coupled spatial fields and channel dynamics.  Guan and
	Queisser~\cite{guan2022calcium} proved well-posedness and analyzed an
	IMEX finite element scheme for reaction--diffusion PDEs coupled to
	interface ODEs.  Gurung and Guan~\cite{gurung2025hybrid} developed a
	hybrid approach in which a neural-network model of channel open
	probabilities is coupled to the calcium PDEs.

	In cardiac electrophysiology, Africa et al.
	combine high-order $Q_p$ spectral elements with BDF1--3 and a linear voltage
	solve~\cite{africa2023matrixfree}.  A recent neuronal tissue
	DG method uses Crank--Nicolson time stepping and includes semi-discrete
	analysis for a simplified ionic problem~\cite{leimersaglio2026dg}.

	Broader electrophysiology research addresses high-order and adaptive
	spatial approximation
	\cite{arthurs2012highorder,cantwell2014spectral,leimersaglio2025padaptive},
	quadrature and ionic interpolation~\cite{krishnamoorthi2013quadrature},
	and time integration, splitting, and multirate or parallel-in-time methods
	\cite{roy2020timestepping,lindner2023splitting,rosilho2024multirate,rosilho2024paralleltime}.
	For FitzHugh--Nagumo systems, Cai et al.~\cite{cai2020fullydiscrete}
	proved first-order-in-time IMEX finite element error bounds.
	Chrysafinos et al.~\cite{chrysafinos2013error} analyzed discontinuous
	Galerkin time stepping with conforming spatial finite elements and
	proved arbitrary-order space--time error estimates for a system of two
	reaction--diffusion PDEs.
	Variable-step BDF2 finite element analysis assumes a bound on adjacent
	step ratios~\cite{liu2024variablebdf2}.  An unconditional
	invariant-set-preserving finite-difference scheme is also available
	\cite{liu2025invariant}.

		In this paper, we construct and analyze constant-step IMEX BDF--finite
	element schemes for the Hodgkin--Huxley model on a branched metric graph.
	Our discretization builds on finite element methods for metric graphs
	\cite{arioli2018finite,weller2025metricgraphs} and IMEX multistep analysis
	for parabolic equations
	\cite{ascher1995imex,akrivis1998imexfem,akrivis2015quasilinear,akrivis2018maxnorm}.
	High-order extrapolation can move gate values outside $[0,1]$ and destroy
	voltage-equation coercivity, while the classical gate rates are not
	globally bounded or Lipschitz in voltage.  We apply a smooth mask to the
	extrapolated gates and a smooth cutoff to the voltage entering the rates.
	The mask preserves a positive total conductance, while the cutoff gives
	bounded Lipschitz rates.  The mask is the identity on $[0,1]$, and the
	cutoff is the identity on a prescribed voltage interval $[-M,M]$
	containing the exact solution range.  The numerical gates remain
	unprojected.
	To our knowledge, unconditional fully discrete energy bounds for
	regularized, linearly implicit BDF schemes of orders three through five
	for the coupled Hodgkin--Huxley cable--gate system on branched axons have
	not previously been established.

	Our main contributions are:
	\begin{enumerate}
	\item We construct linearly implicit BDF$q$--$P_p$ schemes for
	$q=1,\ldots,5$ and any fixed $p\ge1$, with one linear voltage solve
	and independent nodal gate updates per timestep, without nonlinear
	iteration.
	\item Using the gate mask and voltage cutoff, we prove unconditional
	stability in discrete history energies for the complete fully discrete
	BDF1--BDF5 schemes under exact spatial integration.
	\item We prove optimal fully discrete voltage and gate error bounds under
	regularity and starting-history assumptions, and recover the same rate
	with numerical quadrature when its full residual defect is $O(h^{p+1})$.
	\item We assess accuracy and efficiency through manufactured solutions
	and action-potential simulations.
	\end{enumerate}

	The stability proof combines BDF multiplier estimates
	\cite{nevanlinna1981multiplier,akrivis2015stability,akrivis2016multipliers}
	with a positive nodal inner product for the gate equations.  Let
	$\mathcal W=\{m,h,n\}$ denote the gates, $I_h$ the degree-$p$ nodal
	interpolation operator, and $\|\cdot\|_h$ the positive weighted norm of
	gate values at the Gauss--Lobatto nodes.  Writing $\|\cdot\|$ for the
	$L^2(\Gamma)$ norm, the convergence estimate under exact integration and
	the stated regularity and starting-history assumptions takes the form
	\[
	\max_n\left(\|V_h^n-V(t_n)\|
	+\sum_{w\in\mathcal W}\|w_h^n-I_hw(t_n)\|_h\right)
	\le C_{p,q}(T)\bigl(h^{p+1}+\Delta t^q\bigr),
	\]
	without a timestep restriction or an $h$--$\Delta t$ coupling.  The nodal
	gate norm equivalence and interpolation give the same rate for the
	reconstructed gates in $L^2(\Gamma)$.  The voltage update uses previous
	gate values and the gate update uses the new voltage, so their error
	bounds are combined sequentially.

	Manufactured solutions recover the predicted orders through BDF5--$P_4$
	on the Y-shaped graph; a further test checks spatial accuracy on a
	two-junction graph with unequal edge lengths.  Action-potential simulations
	illustrate pulse splitting at a junction, and a controlled work--precision
	study compares members of the proposed family along the tested refinement
	paths.

	The algorithm and exact-integration analysis extend to any fixed
	connected finite metric graph with standard continuity and Kirchhoff
	vertex conditions and uniform cable and gate laws
	(Remark~\ref{rem:general_graph}).  Section~2 presents the model,
	Section~3 defines the schemes and regularizations, Sections~4 and~5
	establish stability and convergence, and Section~6 reports the numerical
	results.  The appendix gives the smooth mask construction.

	\section{Mathematical Model}

	\subsection{Domain and Geometry}
\begin{figure}[htbp]
		\centering
		\begin{tikzpicture}[
			scale=4,
			domain-line/.style={line width=2.5pt, color=blue!60!black},
			stimulus/.style={line width=4pt, color=red!80, opacity=0.7},
			axis/.style={thick, ->, gray!80},
			node-dot/.style={circle, fill=black, inner sep=1.5pt}
			]

			\draw[axis] (-1.2, 0) -- (1.0, 0) node[right, text=black] {$x$};
			\draw[axis] (0, -0.9) -- (0, 0.9) node[above, text=black] {$y$};

			\draw[domain-line] (-1, 0) -- (0, 0)
			node[midway, below=4pt, text=black] {$A_0$};

			\draw[domain-line] (0, 0) -- ({1/sqrt(2)}, {1/sqrt(2)})
			node[midway, above left=2pt, text=black] {$A_1$};

			\draw[domain-line] (0, 0) -- ({1/sqrt(2)}, {-1/sqrt(2)})
			node[midway, below left=2pt, text=black] {$A_2$};

			\draw[stimulus] (-1, 0) -- (-0.9, 0);
			\draw[decorate, decoration={brace, amplitude=4pt, mirror}, thick, red!80!black]
			(-1.0, -0.05) -- (-0.9, -0.05)
			node[midway, below=6pt, text=black] {$I_{\text{ext}}$};

			\node[node-dot] at (-1, 0) {};
			\node[above left=1pt] at (-1, 0) {$(-1,0,0)$};

			\node[node-dot] at (0, 0) {};
			\node[below left=2pt] at (0, 0) {$(0,0,0)$};

			\node[node-dot] at ({1/sqrt(2)}, {1/sqrt(2)}) {};
			\node[right=2pt] at ({1/sqrt(2)}, {1/sqrt(2)}) {$\left(\frac{1}{\sqrt{2}}, \frac{1}{\sqrt{2}},0\right)$};

			\node[node-dot] at ({1/sqrt(2)}, {-1/sqrt(2)}) {};
			\node[right=2pt] at ({1/sqrt(2)}, {-1/sqrt(2)}) {$\left(\frac{1}{\sqrt{2}}, -\frac{1}{\sqrt{2}},0\right)$};

			\draw[-{Latex[length=2mm]}, thick, dashed, orange!90!black]
			(-0.5, 0.0) -- (-0.2, 0.0) node[midway, above=1pt] {\small Wave};

			\draw[-{Latex[length=2mm]}, thick, dashed, orange!90!black]
			(0.1, 0.1) -- (0.35, 0.35);

			\draw[-{Latex[length=2mm]}, thick, dashed, orange!90!black]
			(0.1, -0.1) -- (0.35, -0.35);

		\end{tikzpicture}
		\caption{Planar representative of the branched axon graph $\Gamma=A_0\cup A_1\cup A_2$ embedded in $\R^3$. The external stimulus $I_{\text{ext}}$ is applied near the left endpoint of the trunk.}
		\label{fig:axon_geometry}
	\end{figure}
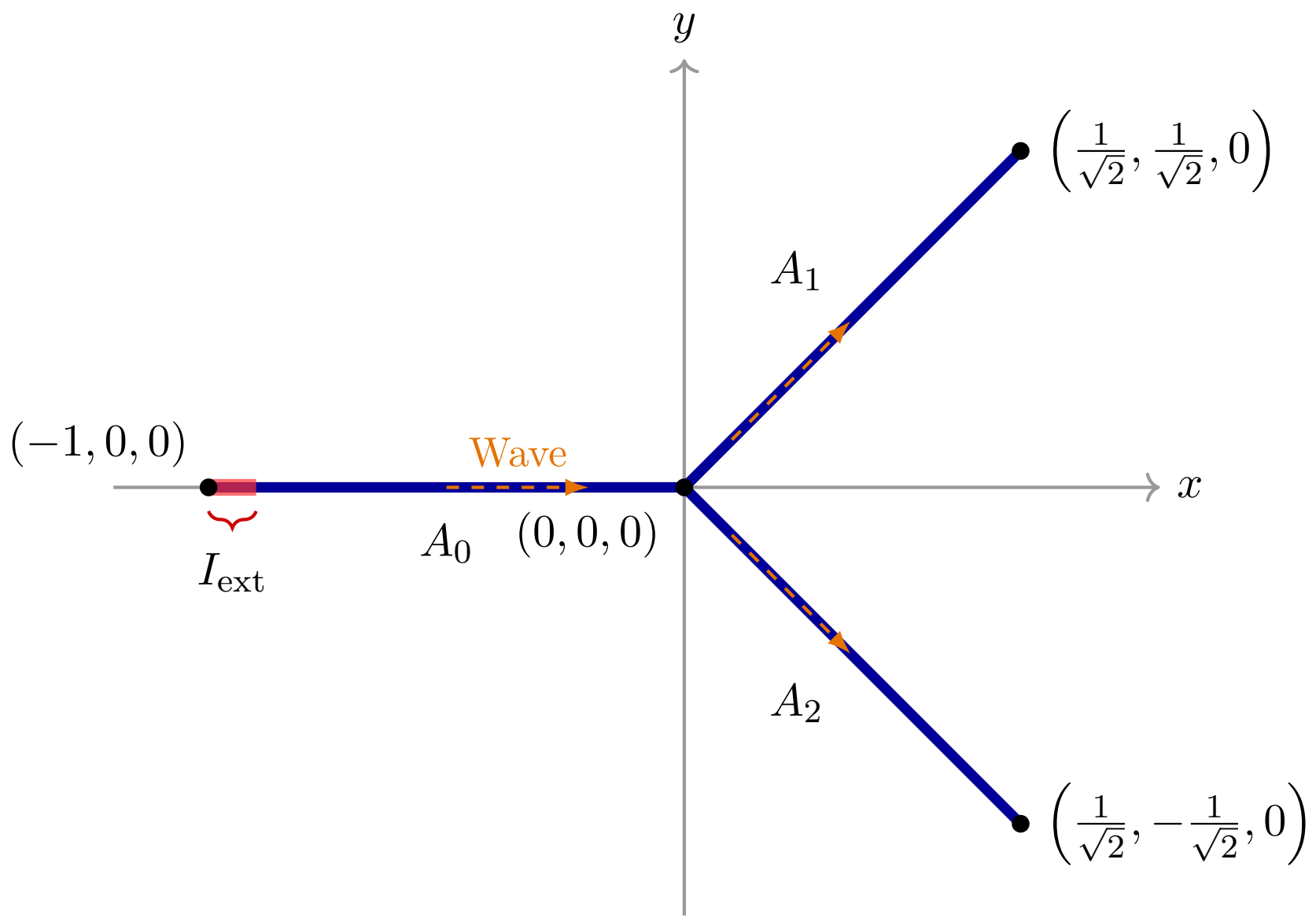

	We use the Y-shaped metric graph
	\[
	\Gamma = A_0 \cup A_1 \cup A_2,
	\]
	with one trunk and two daughter branches, shown in
	Figure~\ref{fig:axon_geometry}.  Its planar embedding in $\R^3$ is
	\begin{align*}
		A_0 &= \{(x,0,0)\;:\; -1\le x\le 0\},\\
		A_1 &= \left\{\left(\frac{s}{\sqrt{2}},\frac{s}{\sqrt{2}},0\right)\;:\; 0\le s\le 1\right\},\\
		A_2 &= \left\{\left(\frac{s}{\sqrt{2}},-\frac{s}{\sqrt{2}},0\right)\;:\; 0\le s\le 1\right\}.
	\end{align*}
	Thus $A_0$ has arc-length coordinate $s=x\in[-1,0]$, while $A_1$ and $A_2$ have arc-length coordinate $s\in[0,1]$. The bifurcation node is located at $(0,0,0)$.
	The edgewise construction also applies to finite metric graphs with
	multiple junctions or cycles, using one shared degree of freedom at each
	vertex.

	\subsection{Cable Equation and Hodgkin--Huxley Dynamics on the Graph}
	Let $V=(V_0,V_1,V_2)$ denote the transmembrane voltage, where each $V_i$ is defined on the edge $A_i$. On each edge, the cable equation reads
	\begin{equation}\label{eq:cable_edge}
		C_m \frac{\partial V_i}{\partial t}
		=
		D \frac{\partial^2 V_i}{\partial s^2}
		-
		I_{\mathrm{ion}}(V_i,m_i,h_i,n_i)
		+
		I_{\mathrm{ext},i}(s,t),
		\qquad \text{on } A_i\times(0,T].
	\end{equation}
	where $D=a/(2R)$.  Throughout the analysis we assume
	\[
	C_m>0,\qquad a>0,\qquad R>0,\qquad
	g_{\mathrm{Na}}\ge0,\qquad g_{\mathrm{K}}\ge0,\qquad g_{\mathrm{L}}>0,
	\]
	so in particular $D>0$.  The ionic current is
	\begin{equation}\label{eq:ion}
		I_{\mathrm{ion}}(V,m,h,n)
		=
		g_{\mathrm{Na}} m^3 h (V-E_{\mathrm{Na}})
		+
		g_{\mathrm{K}} n^4 (V-E_{\mathrm{K}})
		+
		g_{\mathrm{L}} (V-E_{\mathrm{L}}).
	\end{equation}
	Let $\mathcal{W}:=\{m,h,n\}$.  The gating variables $w\in\mathcal{W}$ satisfy
	\begin{equation}\label{eq:gating}
		\frac{\partial w}{\partial t}
		=
		\alpha_w(V)(1-w)-\beta_w(V)w
		=
		\alpha_w(V)-\lambda_w(V)w,
		\qquad \lambda_w(V):=\alpha_w(V)+\beta_w(V).
	\end{equation}
	We use the classical Hodgkin--Huxley rate functions \cite{hodgkin1952quantitative}, in the shifted voltage convention where the resting potential is $V=0$:
	\begin{align*}
		\alpha_m(V)&=\frac{0.1(25-V)}{\exp((25-V)/10)-1},
		&
		\beta_m(V)&=4\exp(-V/18),
		\\
		\alpha_h(V)&=0.07\exp(-V/20),
		&
		\beta_h(V)&=\frac{1}{\exp((30-V)/10)+1},
		\\
		\alpha_n(V)&=\frac{0.01(10-V)}{\exp((10-V)/10)-1},
		&
		\beta_n(V)&=0.125\exp(-V/80).
	\end{align*}
	Standard L'H\^opital regularizations are used at the removable singularities, so that $\alpha_m(25)=1$ and $\alpha_n(10)=0.1$. These rates are nonnegative and hence $\lambda_w(V)>0$ for each gate $w\in\mathcal{W}$.  Moreover, the interval $[0,1]$ is positively invariant for each exact gate equation: at $w=0$ the right-hand side is $\alpha_w(V)\ge0$, while at $w=1$ it is $-\beta_w(V)\le0$.
	At the outer endpoints we impose sealed-end conditions
	\begin{equation}\label{eq:outer_bc}
		\ds V_0(-1,t)=0,\qquad \ds V_1(1,t)=0,\qquad \ds V_2(1,t)=0.
	\end{equation}
	At the junction we impose continuity
	\begin{equation}\label{eq:cont_node}
		V_0(0,t)=V_1(0,t)=V_2(0,t),
	\end{equation}
	and Kirchhoff flux balance
	\begin{equation}\label{eq:kirchhoff}
		D\,\ds V_0(0,t)=D\,\ds V_1(0,t)+D\,\ds V_2(0,t).
	\end{equation}
	Here $\ds$ denotes differentiation in the edge coordinates specified above; with this convention, \eqref{eq:kirchhoff} is the zero-net-flux condition at the junction.

	\subsection{Weak Formulation on the Metric Graph}
	Define
	\[
	H^1(\Gamma)
	:=
	\left\{
	v=(v_0,v_1,v_2)\;:\;
	v_i\in H^1(A_i),\;
	v_0(0)=v_1(0)=v_2(0)
	\right\}.
	\]
	For every integer $k\ge1$, we use the graph Sobolev space
	\[
	H^k(\Gamma)
	:=
	\left\{
	v\in H^1(\Gamma)\;:\;
	v_i\in H^k(A_i),\quad i=0,1,2
	\right\},
	\qquad
	\norm{v}_{H^k(\Gamma)}^2
	:=
	\sum_{i=0}^2\norm{v_i}_{H^k(A_i)}^2.
	\]
	In particular, every element of $H^k(\Gamma)$ has a single common trace at
	the junction, while its higher derivatives are understood edgewise.
	The continuous problem is supplied with initial data
	\begin{equation}\label{eq:continuous_initial_data}
		V(\cdot,0)=V^0\in H^1(\Gamma),
		\qquad
		w(\cdot,0)=w^0\in H^1(\Gamma),
		\quad w\in\mathcal{W},
	\end{equation}
	where $0\le w^0\le1$ on $\Gamma$; compatibility required by the
	regularity assumptions in the convergence analysis is understood.
	For $u,v\in L^2(\Gamma)$, define
	\[
	(u,v)_\Gamma
	:=
	\sum_{i=0}^2 \int_{A_i} u_i(s)v_i(s)\,ds,
	\qquad
	\norm{u}_{L^2(\Gamma)}^2 := (u,u)_\Gamma.
	\]
	For brevity, we write $\norm{\cdot}:=\norm{\cdot}_{L^2(\Gamma)}$ below.\\
	The weak form of \eqref{eq:cable_edge} is: find $V(t)\in H^1(\Gamma)$ such that
	\begin{equation}\label{eq:weak_cont}
		C_m(\partial_t V,v)_\Gamma
		+
		D(\ds V,\ds v)_\Gamma
		+
		(I_{\mathrm{ion}}(V,m,h,n),v)_\Gamma
		=
		(I_{\mathrm{ext}},v)_\Gamma,
		\qquad \forall v\in H^1(\Gamma).
	\end{equation}

	\subsection{External Stimulus}
	For the analysis, $I_{\mathrm{ext}}$ is a prescribed source term satisfying
	\begin{equation}\label{eq:external_stimulus_bound}
		\sup_{0\le t\le T}
		\norm{I_{\mathrm{ext}}(\cdot,t)}
		<\infty .
	\end{equation}
	The compactly supported stimulus used in the computations is specified in
	Section~\ref{sec:numerical_experiments}.

\section{Numerical Method}

Throughout a computation, the independently selected degree $p\ge1$
and BDF order $q\in\{1,\ldots,5\}$ remain fixed, and the time step is
constant.  The finite element
construction is edgewise, with continuity imposed at graph vertices
\cite{arioli2018finite}.

\subsection{Finite Element Space and Nodal Norm}
Let $\mathcal{T}_h$ be a partition of each edge into intervals, compatible
at the graph vertices, and set $h:=\max_{K\in\mathcal{T}_h}|K|$.  Define
\[
X_h^p:=\{v_h\in H^1(\Gamma):v_h|_K\in\mathbb{P}_p(K)
\text{ for every }K\in\mathcal{T}_h\},
\]
where $\mathbb{P}_p(K)$ is the space of polynomials of degree at most $p$.
On the reference interval $[-1,1]$, we use the $p+1$ Gauss--Lobatto nodes
\cite[Eq.~(2.3.12)]{canuto2006spectral}
\[
-1=\xi_0<\xi_1<\cdots<\xi_p=1.
\]
For $p\ge2$, the interior nodes are the zeros of $L_p'$, where $L_p$ is
the Legendre polynomial normalized by $L_p(1)=1$.  The associated weights
are
\[
\omega_i^{(p)}
=\frac{2}{p(p+1)L_p(\xi_i)^2}>0,
\qquad i=0,\ldots,p.
\]
Mapping these nodes affinely to each element gives the global nodal set
$\mathcal{N}_h^p$.  Nodes at a common element endpoint or graph vertex are
identified.  In particular, all incident branches share one junction
degree of freedom.

Let $a_{K,i}$ denote the image of $\xi_i$ on $K$.  The assembled positive
nodal weights are
\[
M_a:=\sum_{\substack{K\in\mathcal{T}_h,\ 0\le i\le p\\a_{K,i}=a}}
\frac{|K|}{2}\omega_i^{(p)},
\qquad a\in\mathcal{N}_h^p.
\]
Thus the junction weight includes contributions from every incident
element.  Since the Gauss--Lobatto rule is exact for constants,
$\sum_{i=0}^p\omega_i^{(p)}=2$ and
$\sum_{a\in\mathcal N_h^p}M_a=|\Gamma|$.
For real nodal vectors $z$ and $y$, define
\begin{equation}\label{eq:lumped_inner_product}
(z,y)_h:=\sum_{a\in\mathcal{N}_h^p}M_a z(a)y(a),
\end{equation}
and
\begin{equation}\label{eq:lumped_nodal_norm}
\norm{z}_h^2:=(z,z)_h.
\end{equation}
We use the same notation for a nodal vector and its unique interpolant in
$X_h^p$.  This inner product will be used for the gate analysis; the
voltage equation below uses the continuous graph inner product.

\begin{lemma}[Equivalence of the nodal and continuous norms]\label{lem:nodal_norm}
For every integer $p\ge1$,
\begin{equation}\label{eq:nodal_norm_equivalence}
\norm{z_h}^2
\le\norm{z_h}_h^2
\le\left(2+\frac1p\right)\norm{z_h}^2,
\qquad z_h\in X_h^p.
\end{equation}
In particular, the equivalence constants can be chosen as $1$ and $3$,
independently of $h$ and $p$.
\end{lemma}

\begin{proof}
On the reference interval, write
$\widehat z=\sum_{j=0}^p\widehat z_jL_j$.
Legendre orthogonality and its Gauss--Lobatto discrete counterpart give
\cite{teukolsky2015mass}
\begin{align*}
\int_{-1}^1|\widehat z(\xi)|^2\,d\xi
&=\sum_{j=0}^p\frac{2}{2j+1}|\widehat z_j|^2,\\
\sum_{i=0}^p\omega_i^{(p)}|\widehat z(\xi_i)|^2
&=\sum_{j=0}^{p-1}\frac{2}{2j+1}|\widehat z_j|^2
+\frac2p|\widehat z_p|^2.
\end{align*}
Only the highest-mode weight differs, by the factor
$(2p+1)/p=2+1/p\le3$.  Affine scaling to each element and assembly
over the graph therefore prove \eqref{eq:nodal_norm_equivalence}.
\end{proof}

\subsection{BDF Time Discretization}
Let $\dt>0$ be constant and $t_n=n\dt$.  Following the BDF formulation
in \cite{akrivis2015stability}, define the newest-first BDF
coefficients $a_j^{(q)}$ by
\begin{equation}\label{eq:bdf_coefficients}
\sum_{j=0}^q a_j^{(q)}\zeta^j
=\sum_{\ell=1}^q\frac{(1-\zeta)^\ell}{\ell}.
\end{equation}
For a sequence $\{z^n\}$, the undivided BDF difference is
\begin{equation}\label{eq:bdf_difference}
\bdf{z}^{\,n+1}:=\sum_{j=0}^q a_j^{(q)}z^{n+1-j},
\qquad n\ge q-1.
\end{equation}
The approximation to the time derivative is
$\bdf{z}^{\,n+1}/\dt$.  In particular,
$a_0^{(q)}=\sum_{\ell=1}^q1/\ell>0$.
The coefficients for the permitted orders are listed in
Table~\ref{tab:bdf_coefficients}.

\begin{table}[htbp]
\centering
\caption{Newest-first coefficients of the undivided BDF difference.}
\label{tab:bdf_coefficients}
\renewcommand{\arraystretch}{1.5}
\begin{tabular}{@{}c r r r r r r@{}}
\toprule
$q$ & $a_0^{(q)}$ & $a_1^{(q)}$ & $a_2^{(q)}$
& $a_3^{(q)}$ & $a_4^{(q)}$ & $a_5^{(q)}$ \\
\midrule
$1$ & $1$ & $-1$ & --- & --- & --- & --- \\
$2$ & $\frac32$ & $-2$ & $\frac12$ & --- & --- & --- \\
$3$ & $\frac{11}{6}$ & $-3$ & $\frac32$ & $-\frac13$ & --- & --- \\
$4$ & $\frac{25}{12}$ & $-4$ & $3$ & $-\frac43$ & $\frac14$ & --- \\
$5$ & $\frac{137}{60}$ & $-5$ & $5$ & $-\frac{10}{3}$ & $\frac54$ & $-\frac15$ \\
\bottomrule
\end{tabular}
\end{table}

\subsection{Gate Extrapolation and Mask}
For each gate $w\in\mathcal{W}$, let $w_h^n\in X_h^p$ denote the
numerical solution at $t_n$.  For $n\ge q-1$, define the order-$q$
extrapolation using the coefficients of the IMEX BDF polynomial in
\cite{akrivis2015stability}:
\begin{equation}\label{eq:extrap}
w_*^{\,n+1}
=\sum_{j=0}^{q-1}b_j^{(q)}w_h^{n-j},
\qquad
b_j^{(q)}:=(-1)^j\binom{q}{j+1}.
\end{equation}
For example, $q=3$ gives
$w_*^{\,n+1}=3w_h^n-3w_h^{n-1}+w_h^{n-2}$.
The extrapolated finite element gates need not remain in $[0,1]$.
To retain coercivity of the voltage equation, we use the fixed $C^6$
mask constructed in Appendix~\ref{app:masking_function}.  It satisfies
\[
\phi(x)=x\quad (0\le x\le1),
\qquad -\varepsilon\le\phi(x)\le1+\varepsilon,
\qquad |\phi(x)-\phi(y)|\le L_\phi|x-y|,
\qquad L_\phi=1.
\]
Choose $\varepsilon>0$ so that
\[
g_{\mathrm{Na}}\varepsilon(1+\varepsilon)^3
\le\frac{g_{\mathrm{L}}}{2}.
\]
These mask properties suffice for every permitted $(q,p)$.

Write
\begin{equation}\label{def_P}
\begin{aligned}
&\bm{\mathcal{P}}=(m,h,n),
\qquad \bm{\mathcal{P}}^n=(m(t_n),h(t_n),n(t_n)),\\
&\bm{\mathcal{P}}_h^n=(m_h^n,h_h^n,n_h^n),
\qquad \bm{\mathcal{P}}_*^{n+1}
=(m_*^{\,n+1},h_*^{\,n+1},n_*^{\,n+1}),\\
&\Phi(\mathbf r):=(\phi(r_m),\phi(r_h),\phi(r_n)),
\qquad\mathbf r=(r_m,r_h,r_n).
\end{aligned}
\end{equation}
For any triple $\mathbf r=(r_m,r_h,r_n)$, define
\begin{equation}\label{def_GF}
\begin{aligned}
G(\mathbf r)
&=g_{\mathrm{Na}}r_m^3r_h+g_{\mathrm{K}}r_n^4+g_{\mathrm{L}},\\
F(\mathbf r)
&=g_{\mathrm{Na}}r_m^3r_hE_{\mathrm{Na}}
+g_{\mathrm{K}}r_n^4E_{\mathrm{K}}+g_{\mathrm{L}}E_{\mathrm{L}}.
\end{aligned}
\end{equation}
These maps are evaluated pointwise on the masked extrapolated finite
element fields $\Phi(\bm{\mathcal{P}}_*^{n+1})$.

\subsection{Voltage Solve}
For $n\ge q-1$, given the voltage and gate histories at levels
$n,n-1,\ldots,n-q+1$, define
\begin{equation}\label{eq:discrete_coefficients}
\begin{aligned}
G_h^{n+1}&:=G(\Phi(\bm{\mathcal P}_*^{n+1})),\\
F_h^{n+1}&:=F(\Phi(\bm{\mathcal P}_*^{n+1})),\\
f_h^{n+1}&:=F_h^{n+1}+I_{\mathrm{ext}}^{n+1},
\end{aligned}
\end{equation}
where $I_{\mathrm{ext}}^{n+1}=I_{\mathrm{ext}}(\cdot,t_{n+1})$.
Find $V_h^{n+1}\in X_h^p$ such that
\begin{equation}\label{eq:voltage_scheme}
\begin{split}
\frac{C_m}{\dt}\bigl(\bdf{V}_h^{\,n+1},\chi_h\bigr)_\Gamma
+D(\ds V_h^{n+1},\ds\chi_h)_\Gamma
&+\bigl(G_h^{n+1}V_h^{n+1},\chi_h\bigr)_\Gamma\\
&=\bigl(f_h^{n+1},\chi_h\bigr)_\Gamma,
\qquad\chi_h\in X_h^p,
\end{split}
\end{equation}
All coefficients involving the gates are known from previous levels, so
this equation is linear in $V_h^{n+1}$.

\subsection{Voltage Cutoff for the Rate Functions}
Let $M>0$ be chosen larger than the $L^\infty$-range of the exact voltage.
Use a $C^6$ cutoff $\psi_M:\R\to\R$ satisfying
\[
\psi_M(x)=x\quad\text{for }|x|\le M,
\qquad |\psi_M(x)|\le M+1,
\qquad |\psi_M'(x)|\le1.
\]
The smooth saturation in Appendix~\ref{app:masking_function}, applied
symmetrically outside $[-M,M]$ with unit saturation width, provides such
a cutoff.  Using the removable-singularity regularizations of the rate
functions, define
\[
\widetilde{\alpha}_w(V):=\alpha_w(\psi_M(V)),
\qquad
\widetilde{\beta}_w(V):=\beta_w(\psi_M(V)),
\qquad
\widetilde{\lambda}_w(V)
:=\widetilde{\alpha}_w(V)+\widetilde{\beta}_w(V).
\]
These functions are globally bounded and Lipschitz, and
$\widetilde{\lambda}_w$ is strictly positive.  If the extrapolated gate
fields remain in $[0,1]$ and the voltage arguments remain in $[-M,M]$,
the mask and cutoff are inactive and the original Hodgkin--Huxley
coefficients and rates are recovered.

\subsection{Nodal Gate Update}
After computing $V_h^{n+1}$, update each $w\in\mathcal{W}$ at every
$a\in\mathcal{N}_h^p$ by
\begin{equation}\label{eq:gate_scheme}
\frac{1}{\dt}\bdf{w}_h^{\,n+1}(a)
=\widetilde{\alpha}_w\bigl(V_h^{n+1}(a)\bigr)
-\widetilde{\lambda}_w\bigl(V_h^{n+1}(a)\bigr)w_h^{n+1}(a),
\qquad n\ge q-1.
\end{equation}
Equivalently,
\begin{equation}\label{eq:gate_unprojected}
w_h^{n+1}(a)
=\frac{
-\displaystyle\sum_{j=1}^q a_j^{(q)}w_h^{n+1-j}(a)
+\dt\,\widetilde{\alpha}_w\bigl(V_h^{n+1}(a)\bigr)
}{
a_0^{(q)}+\dt\,\widetilde{\lambda}_w\bigl(V_h^{n+1}(a)\bigr)
}.
\end{equation}
The denominator is positive for every $\dt>0$.  These nodal values,
including the single value at the shared junction node, therefore
determine a unique $w_h^{n+1}\in X_h^p$.
Since $V_h^{n+1}$ is already known, \eqref{eq:gate_unprojected} is
evaluated directly.  Thus the complete timestep consists of one global
linear voltage solve and independent scalar gate updates, with no
nonlinear iteration.

\begin{remark}
The gate update does not project the numerical values onto $[0,1]$.
\end{remark}

\subsection{Initialization}
Prescribe $V_h^j,w_h^j\in X_h^p$ for $j=0,\ldots,q-1$ and
$w\in\mathcal W$.  The convergence analysis requires
Assumption~\ref{ass:starting_errors}; the localized-stimulus experiment
uses the exact resting histories \eqref{eq:initialization}.

\section{Stability Analysis}

We establish unconditional discrete-energy stability for the complete
voltage--gate scheme, simultaneously for the BDF orders
$q\in\{1,\ldots,5\}$ and the conforming finite element degrees
$p\ge1$.  Throughout this section the inner products in
\eqref{eq:voltage_scheme} are evaluated exactly.  Neither estimate requires
a numerical maximum principle or a restriction on the time step; dependence
on the starting values is expressed through discrete history energies.

\subsection{Basic Bounds Induced by the Mask}
The pointwise mask bounds imply
\[
\phi(m_*^{n+1})^3\phi(h_*^{n+1})
\ge-\varepsilon(1+\varepsilon)^3.
\]
And $\phi(n_*^{n+1})^4$ is nonnegative.  Therefore the prescribed
choice $g_{\mathrm{Na}}\varepsilon(1+\varepsilon)^3\le g_{\mathrm L}/2$
gives
\begin{equation}\label{eq:G_bounds}
\begin{aligned}
0<g_0:=\frac{g_{\mathrm L}}2
&\le G_h^{n+1}\le g_1,\\
g_1&:=g_{\mathrm L}
 +(g_{\mathrm{Na}}+g_{\mathrm K})(1+\varepsilon)^4,
\qquad n\ge q-1.
\end{aligned}
\end{equation}
$G_h^{n+1}$ is defined in \eqref{eq:discrete_coefficients}. Similarly, from \eqref{def_GF} and \eqref{eq:discrete_coefficients}, we have
\[
\abs{F_h^{n+1}}
\le g_{\mathrm{Na}}(1+\varepsilon)^4\abs{E_{\mathrm{Na}}}
 +g_{\mathrm K}(1+\varepsilon)^4\abs{E_{\mathrm K}}
 +g_{\mathrm L}\abs{E_{\mathrm L}}=:C_E.
\]
Consequently,
\begin{equation}\label{eq:F_bound}
\norm{F_h^{n+1}}
\le |\Gamma|^{1/2}C_E.
\end{equation}
Together with \eqref{eq:external_stimulus_bound}, this yields a constant
$C_F$ independent of $h$, $\dt$, $q$, $p$, and the numerical gate values,
such that
\begin{equation}\label{eq:f_bound}
\norm{f_h^{n+1}}\le C_F,
\qquad n\ge q-1,\quad t_{n+1}\le T.
\end{equation}
In particular, these bounds require no control of the unmasked gate
extrapolations.

\subsection{BDF Multiplier Inequality}
We use the admissible Nevanlinna--Odeh multipliers
\begin{equation}\label{eq:multiplier_values}
(\eta_1,\eta_2,\eta_3,\eta_4,\eta_5)
=(0,0,0.0836,0.2879,0.8160).
\end{equation}
The multiplier construction follows from Dahlquist's $G$-stability theory
and the Nevanlinna--Odeh criterion
\cite{dahlquist1978g,nevanlinna1981multiplier,akrivis2015stability}.
Precise values of the optimal multipliers are given in
\cite{akrivis2016multipliers}.  In particular, the optimal fourth-order
value is approximately $0.287806557$, so the exact decimal $0.2879$ in
\eqref{eq:multiplier_values} is chosen above that threshold.

\begin{lemma}[BDF multiplier inequality]\label{lem:bdf_multiplier}
For each $q\in\{1,\ldots,5\}$, there is a symmetric positive-definite
matrix $\mathbf M_q\in\R^{q\times q}$ such that every sequence in a real
inner product space $X$ satisfies
\begin{equation}\label{eq:NO_multiplier}
\bigl(\bdf{z}^{\,n+1},z^{n+1}-\eta_qz^n\bigr)_X
\ge\frac12\left(
 \norm{\bm{\mathcal Z}_q^{n+1}}_{\mathbf M_q,X}^2
 -\norm{\bm{\mathcal Z}_q^n}_{\mathbf M_q,X}^2
\right),
\qquad n\ge q-1,
\end{equation}
where
\[
\bm{\mathcal Z}_q^n:=[z^n,z^{n-1},\ldots,z^{n-q+1}]^T,
\qquad
\norm{\bm{\mathcal Z}_q^n}_{\mathbf M_q,X}^2
:=\sum_{i,j=1}^q(M_q)_{ij}
 (z^{n+1-i},z^{n+1-j})_X.
\]
Moreover, there are constants $c_*>0$ and $c^*<\infty$, independent of
$q\in\{1,\ldots,5\}$ and of $X$, such that
\begin{equation}\label{eq:norm_equiv}
c_*\sum_{j=0}^{q-1}\norm{z^{n-j}}_X^2
\le\norm{\bm{\mathcal Z}_q^n}_{\mathbf M_q,X}^2
\le c^*\sum_{j=0}^{q-1}\norm{z^{n-j}}_X^2.
\end{equation}
\end{lemma}

\begin{proof}
Apply the Nevanlinna--Odeh identity recalled by
Akrivis~\cite[Lemma~1.1]{akrivis2015stability} and discard its
nonnegative remainder.  For the admissible multipliers in
\eqref{eq:multiplier_values}, let $\mathbf G_q$ denote the positive-definite
matrix in that identity and let $J_q$ reverse the history order.  Taking
$\mathbf M_q=2J_q\mathbf G_qJ_q$ gives our newest-first convention and
the factor $1/2$ in \eqref{eq:NO_multiplier}.  The norm equivalence
\eqref{eq:norm_equiv} follows from positive definiteness; its constants
can be chosen independently of $X$ and uniformly over the five orders
$q=1,\ldots,5$.
\end{proof}

We omit the subscript $X$ when the underlying inner product is
$(\cdot,\cdot)_\Gamma$.

\subsection{Unconditional Voltage Stability}
\begin{theorem}[Masked-coefficient voltage stability]\label{thm:stability}
Fix $q\in\{1,\ldots,5\}$ and an integer $p\ge1$.  Let
$\{V_h^n\}\subset X_h^p$ satisfy \eqref{eq:voltage_scheme}, with starting
values $V_h^0,\ldots,V_h^{q-1}\in X_h^p$.  Define
\begin{equation}\label{eq:energy_def}
\begin{aligned}
\bm{\mathcal V}_{h,q}^n
&:=[V_h^n,V_h^{n-1},\ldots,V_h^{n-q+1}]^T,\\
\mathcal E_V^n
&:=\norm{\bm{\mathcal V}_{h,q}^n}_{\mathbf M_q}^2
 +\frac{D\eta_q^2\dt}{C_m}\norm{\ds V_h^n}^2,
\qquad n\ge q-1.
\end{aligned}
\end{equation}
There exists $C(T)>0$, independent of $h$, $\dt$, $p$, and $q$, such
that, for every $N\ge q-1$ with
$t_N\le T$,
\begin{equation}\label{eq:stability_result}
\max_{q-1\le m\le N}\norm{V_h^m}^2
\le C(T)\left(\mathcal E_V^{q-1}+TC_F^2\right).
\end{equation}
The same right-hand side, with a possibly enlarged $C(T)$, bounds
$\max\limits_{q-1\le m\le N}\mathcal E_V^m$.
\end{theorem}

\begin{proof}
For $n\ge q-1$, put
\[
\chi_h^{n+1}:=V_h^{n+1}-\eta_qV_h^n\in X_h^p.
\]
Testing \eqref{eq:voltage_scheme} with $\chi_h^{n+1}$ gives
\begin{equation}\label{eq:test_eq}
\begin{aligned}
\frac{C_m}{\dt}
 (\bdf{V}_h^{\,n+1},\chi_h^{n+1})_\Gamma
 &+D(\ds V_h^{n+1},\ds\chi_h^{n+1})_\Gamma\\
 &+(G_h^{n+1}V_h^{n+1},\chi_h^{n+1})_\Gamma
 =(f_h^{n+1},\chi_h^{n+1})_\Gamma.
\end{aligned}
\end{equation}

\medskip\noindent{Step 1: the temporal and diffusion terms.}
Lemma~\ref{lem:bdf_multiplier} yields
\begin{equation}\label{eq:temporal_est}
(\bdf{V}_h^{\,n+1},\chi_h^{n+1})_\Gamma
\ge\frac12\left(
\norm{\bm{\mathcal V}_{h,q}^{n+1}}_{\mathbf M_q}^2
-\norm{\bm{\mathcal V}_{h,q}^{n}}_{\mathbf M_q}^2\right).
\end{equation}
The identity $2a(a-\eta_qb)=a^2-\eta_q^2b^2+(a-\eta_qb)^2$ gives
\begin{equation}\label{eq:diff_est}
D(\ds V_h^{n+1},\ds\chi_h^{n+1})_\Gamma
=\frac D2\left(
\norm{\ds V_h^{n+1}}^2-\eta_q^2\norm{\ds V_h^n}^2
+\norm{\ds\chi_h^{n+1}}^2\right).
\end{equation}

\medskip\noindent{Step 2: the reaction and source terms.}
The same identity, applied pointwise with the positive weight
$G_h^{n+1}$, gives
\begin{align}
(G_h^{n+1}V_h^{n+1},\chi_h^{n+1})_\Gamma
&=\frac12\Bigl(
 \norm{\sqrt{G_h^{n+1}}V_h^{n+1}}^2
 -\eta_q^2\norm{\sqrt{G_h^{n+1}}V_h^n}^2 
 +\norm{\sqrt{G_h^{n+1}}\chi_h^{n+1}}^2\Bigr)\notag\\
&\ge\frac{g_0}2\norm{V_h^{n+1}}^2
 -\frac{\eta_q^2g_1}2\norm{V_h^n}^2
 +\frac{g_0}2\norm{\chi_h^{n+1}}^2.
\label{eq:react_est}
\end{align}
By Cauchy--Schwarz and Young's inequality,
\begin{equation}\label{eq:rhs_est}
(f_h^{n+1},\chi_h^{n+1})_\Gamma
\le\frac1{2g_0}\norm{f_h^{n+1}}^2
 +\frac{g_0}2\norm{\chi_h^{n+1}}^2.
\end{equation}
The last term is absorbed by \eqref{eq:react_est}.

\medskip\noindent{Step 3: the energy recurrence.}
Combining \eqref{eq:test_eq}--\eqref{eq:rhs_est}, multiplying by
$2\dt/C_m$, and dropping the nonnegative term involving
$\norm{\ds\chi_h^{n+1}}^2$ yields
\begin{align}
&\norm{\bm{\mathcal V}_{h,q}^{n+1}}_{\mathbf M_q}^2
 -\norm{\bm{\mathcal V}_{h,q}^n}_{\mathbf M_q}^2
 +\frac{D\dt}{C_m}
 \left(\norm{\ds V_h^{n+1}}^2
       -\eta_q^2\norm{\ds V_h^n}^2\right)
 +\frac{g_0\dt}{C_m}\norm{V_h^{n+1}}^2\notag\\
&\hspace{2cm}\le\frac{\dt}{C_mg_0}\norm{f_h^{n+1}}^2
 +\frac{\eta_q^2g_1\dt}{C_m}\norm{V_h^n}^2.
\label{eq:combined_est}
\end{align}
By \eqref{eq:norm_equiv},
\[
\norm{V_h^n}^2
\le c_*^{-1}\norm{\bm{\mathcal V}_{h,q}^n}_{\mathbf M_q}^2
\le c_*^{-1}\mathcal E_V^n.
\]
The definition \eqref{eq:energy_def} therefore transforms
\eqref{eq:combined_est} into
\begin{equation}\label{eq:energy_rec}
\begin{aligned}
\mathcal E_V^{n+1}
 &+\frac{D(1-\eta_q^2)\dt}{C_m}\norm{\ds V_h^{n+1}}^2
 +\frac{g_0\dt}{C_m}\norm{V_h^{n+1}}^2
\\
&\le(1+C_v\dt)\mathcal E_V^n
 +\frac{\dt}{C_mg_0}\norm{f_h^{n+1}}^2,
\end{aligned}
\end{equation}
where
\[
C_v:=\frac{\eta_q^2g_1}{C_mc_*}.
\]
For every order considered, $0\le\eta_q<1$, so both dissipation terms on
the left are nonnegative.  For $q=1,2$, $\eta_q=0$ and
$C_v=0$; the history energy alone is sufficient, and the
diffusion term remains nonnegative.

\medskip\noindent{Step 4: iteration.}
Dropping the nonnegative dissipation terms in \eqref{eq:energy_rec},
iterating from $n=q-1$, and using
$1+C_v\dt\le\exp(C_v\dt)$ and \eqref{eq:f_bound} gives,
for $q-1\le m\le N$,
\begin{align}
\mathcal E_V^m
&\le\exp\bigl(C_v(t_m-t_{q-1})\bigr)
 \left(\mathcal E_V^{q-1}
 +\frac1{C_mg_0}\sum_{n=q-1}^{m-1}\dt\norm{f_h^{n+1}}^2\right)\notag\\
&\le\exp(C_vT)
 \left(\mathcal E_V^{q-1}+\frac{TC_F^2}{C_mg_0}\right).
\label{eq:gronwall}
\end{align}
The empty sum for $m=q-1$ is zero.  Since
$c_*\norm{V_h^m}^2\le\mathcal E_V^m$, this proves
\eqref{eq:stability_result}. 
\end{proof}

\begin{remark}
The estimate is unconditional with respect to $h$ and $\dt$.  For
$q=3,4,5$, its starting energy contains
$D\eta_q^2\dt\norm{\ds V_h^{q-1}}^2/C_m$, so it is a discrete-energy
estimate and does not assert a pure $L^2$-to-$L^2$ bound for arbitrary
starting histories.  For $q=1,2$ that additional term vanishes.
Under the resting initialization \eqref{eq:initialization}, all voltage
starting values are zero and $\mathcal E_V^{q-1}=0$; hence
\eqref{eq:stability_result} reduces to
$\max_{q-1\le m\le N}\norm{V_h^m}^2\le C(T)TC_F^2$.
\end{remark}

\subsection{Unconditional Gate Stability}

\begin{theorem}[Gate stability]\label{thm:gate_stability}
Fix $q\in\{1,\ldots,5\}$ and an integer $p\ge1$.  Let the voltage
sequence be arbitrary, and let each gate sequence satisfy
\eqref{eq:gate_scheme} with prescribed starting values
$w_h^0,\ldots,w_h^{q-1}$.  For $n\ge q-1$, define
\[
\bm{\mathcal W}_{h,q}^n
:=[w_h^n,w_h^{n-1},\ldots,w_h^{n-q+1}]^T,
\qquad
\mathcal E_w^n
:=\norm{\bm{\mathcal W}_{h,q}^n}_{\mathbf M_q,h}^2,
\]
where the history norm is induced by $(\cdot,\cdot)_h$.  There is a
constant $C(T)>0$, independent of $h$, $\dt$, $p$, and $q$, such that
\begin{equation}\label{eq:gate_stability_result}
\max_{q-1\le m\le N}\sum_{w\in\mathcal W}\mathcal E_w^m
\le C(T)\left(
\sum_{w\in\mathcal W}\mathcal E_w^{q-1}+T\right)
\end{equation}
whenever $N\ge q-1$ and $t_N\le T$.  Consequently,
$\max_{q-1\le m\le N}\sum_{w\in\mathcal W}\norm{w_h^m}_h^2$ and the
corresponding reconstructed $L^2(\Gamma)$ norms obey the same bound.
\end{theorem}

\begin{proof}
The cutoff construction gives constant bounds, common to the three gates, such
that
\begin{equation}\label{eq:regularized_rate_bounds}
0\le\widetilde\alpha_w(v)\le\alpha_{1},
\qquad
0<\lambda_{0}
\le\widetilde\lambda_w(v)
\le\lambda_{1},
\qquad v\in\mathbb R.
\end{equation}
Fix $w\in\mathcal W$ and set
\[
\zeta_h^{n+1}:=w_h^{n+1}-\eta_qw_h^n,
\quad
\alpha_{w,h}^{n+1}:=\widetilde\alpha_w(V_h^{n+1}),
\quad
\lambda_{w,h}^{n+1}:=\widetilde\lambda_w(V_h^{n+1}).
\]
Testing \eqref{eq:gate_scheme} with $\zeta_h^{n+1}$ in
$(\cdot,\cdot)_h$, as Theorem \ref{thm:stability}, Lemma~\ref{lem:bdf_multiplier} and the weighted square
identity give
\begin{align*}
\frac{1}{2\dt}\left(\mathcal E_w^{n+1}-\mathcal E_w^n\right)
&+\frac{\lambda_{0}}2\norm{w_h^{n+1}}_h^2
-\frac{\eta_q^2\lambda_{1}}2\norm{w_h^n}_h^2\\
&+\frac{\lambda_{0}}2\norm{\zeta_h^{n+1}}_h^2
\le(\alpha_{w,h}^{n+1},\zeta_h^{n+1})_h.
\end{align*}
Young's inequality and $\sum_{a\in\mathcal N_h^p}M_a=|\Gamma|$ yield
\[
(\alpha_{w,h}^{n+1},\zeta_h^{n+1})_h
\le\frac{|\Gamma|\alpha_{1}^2}
{2\lambda_{0}}
+\frac{\lambda_{0}}2
\norm{\zeta_h^{n+1}}_h^2.
\]
After cancellation and use of \eqref{eq:norm_equiv},
\begin{equation}\label{eq:gate_energy_rec}
\mathcal E_w^{n+1}
\le(1+C_{g}\dt)\mathcal E_w^n
+\frac{|\Gamma|\alpha_{1}^2}
{\lambda_{0}}\dt,
\qquad
C_{g}
:=\frac{\eta_q^2\lambda_{1}}{c_*}.
\end{equation}
Iterating \eqref{eq:gate_energy_rec} from $n=q-1$, using
$1+C_{g}\dt\le\exp(C_{g}\dt)$, and summing
over the three gates proves \eqref{eq:gate_stability_result}.  The final assertions
follow from \eqref{eq:norm_equiv} and Lemma~\ref{lem:nodal_norm}.
Uniformity over $q$ and $p$ follows from the fixed multiplier and
nodal-norm constants.
\end{proof}

Combining Theorems~\ref{thm:stability} and \ref{thm:gate_stability} gives
unconditional stability of the complete voltage--gate scheme:
\[
\max_{q-1\le m\le N}
\left(\mathcal E_V^m+
\sum_{w\in\mathcal W}\mathcal E_w^m\right)
\le C(T)\left(
\mathcal E_V^{q-1}
+\sum_{w\in\mathcal W}\mathcal E_w^{q-1}
+T(1+C_F^2)\right).
\]

\section{Convergence Analysis}
\label{sec:convergence}

Fix a time order $q\in\{1,\ldots,5\}$ and an integer polynomial degree
$p\ge1$.  All time steps in this section are constant.
We write $I_h$ for the degree-$p$ nodal interpolation operator into
$X_h^p$ and use the nodal norm \eqref{eq:lumped_nodal_norm} for gate
errors, with the norm equivalence from Lemma~\ref{lem:nodal_norm}.

For the exact gate $w\in\mathcal W$, define
\begin{equation}\label{eq:gate_error_definition}
w_I^n:=I_hw(t_n),
\qquad
e_w^n:=w_h^n-w_I^n.
\end{equation}
Unless a different dependence is stated, constants below may depend on
the graph, the physical parameters, the fixed regularizations, the
regularity bounds, and $p,q$, but not on $h$ or $\dt$.

\begin{assumption}[Regularity and boundedness]
\label{ass:conv_regularity}
The exact solution satisfies
\[
V,w\in W^{q+1,\infty}(0,T;H^{p+1}(\Gamma)),
\qquad w\in\mathcal W.
\]
The graph Sobolev norms are understood edgewise, with the common
junction trace prescribed by $H^1(\Gamma)$.  The exact initial gates
satisfy $0\le w(\cdot,0)\le1$.  Moreover, there is a constant $M_0<M$
such that
\[
|V(x,t)|\le M_0,
\qquad (x,t)\in\Gamma\times[0,T].
\]
\end{assumption}
The gate mask acts as the identity on the exact
gates, and the voltage cutoff acts as the identity on the exact voltage.
The regularized rates satisfy \eqref{eq:regularized_rate_bounds} and,
for a constant $L_M$ common to the three gates,
\begin{equation}\label{eq:rate_lipschitz}
|\widetilde\alpha_w(u)-\widetilde\alpha_w(v)|
+|\widetilde\lambda_w(u)-\widetilde\lambda_w(v)|
\le L_M|u-v|,
\qquad u,v\in\mathbb R.
\end{equation}

\subsection{Interpolation, extrapolation, and the Ritz projection}

The standard degree-$p$ interpolation estimate on each edge \cite{thomee2006galerkin}, summed over
the graph, is 
\begin{equation}\label{eq:graph_interpolation_error}
\norm{v-I_hv}+h\norm{\ds(v-I_hv)}
\le Ch^{p+1}\norm{v}_{H^{p+1}(\Gamma)}.
\end{equation}
The common junction trace ensures that $I_hv\in X_h^p$.

\begin{lemma}[Extrapolation and coefficient consistency]
\label{lem:projection_extrapolation}
Under Assumption~\ref{ass:conv_regularity}, for each $w\in\mathcal W$ and
$n\ge q-1$ with $t_{n+1}\le T$, the extrapolated gate satisfies
\begin{equation}\label{eq:extrap_error}
\norm{w_*^{n+1}-w_I^{n+1}}_h
\le C\sum_{j=0}^{q-1}\norm{e_w^{n-j}}_h+C\dt^q.
\end{equation}
Here $e_w^n$ is defined in \eqref{eq:gate_error_definition}.
Let the exact coefficients at $t_{n+1}$ be
\begin{equation}\label{eq:exact_coefficients}
G^{n+1}:=G(\bm{\mathcal P}^{n+1}),
\qquad F^{n+1}:=F(\bm{\mathcal P}^{n+1}).
\end{equation}
With $G_h^{n+1}$ and $F_h^{n+1}$ defined in
\eqref{eq:discrete_coefficients}, we also have
\begin{equation}\label{eq:coefficient_error}
\norm{G_h^{n+1}-G^{n+1}}+\norm{F_h^{n+1}-F^{n+1}}
\le C\sum_{w\in\mathcal W}\sum_{j=0}^{q-1}
\norm{e_w^{n-j}}_h+C(h^{p+1}+\dt^q).
\end{equation}
\end{lemma}

\begin{proof}
Recall $b_j^{(q)}=(-1)^j\binom q{j+1}$ for $0\le j\le q-1$.
Polynomial extrapolation from the previous $q$ time levels is exact for
polynomials of degree at most $q-1$.  Let $C^0(\Gamma)$ denote the space of
continuous functions on $\Gamma$ with the supremum norm.  Taylor expansion at
$t_{n+1}$ gives, at every node $a\in\mathcal N_h^p$,
\[
\left|\sum_{j=0}^{q-1} b_j^{(q)}w(a,t_{n-j})-w(a,t_{n+1})\right|
\le C\dt^q
\norm{\partial_t^q w}_{L^\infty(t_{n+1-q},t_{n+1};C^0(\Gamma))}.
\]
The right-hand side is bounded by $C\dt^q$ under
Assumption~\ref{ass:conv_regularity}, since the one-dimensional Sobolev
embedding applies on every edge.  The positive nodal weights sum to
$|\Gamma|$, so the same estimate holds in $\norm{\cdot}_h$.
The identity
\[
w_*^{n+1}-w_I^{n+1}
=\sum_{j=0}^{q-1} b_j^{(q)}e_w^{n-j}
+\sum_{j=0}^{q-1} b_j^{(q)}w_I^{n-j}-w_I^{n+1}
\]
therefore proves \eqref{eq:extrap_error}.

Since $\phi(w(t_{n+1}))=w(t_{n+1})$, global Lipschitz continuity of $\phi$,
Lemma~\ref{lem:nodal_norm}, and \eqref{eq:graph_interpolation_error} give
\begin{align}
\norm{\phi(w_*^{n+1})-w(t_{n+1})}
&\le L_\phi\norm{w_*^{n+1}-w(t_{n+1})}\notag\\
&\le C\norm{w_*^{n+1}-w_I^{n+1}}_h
+C\norm{w_I^{n+1}-w(t_{n+1})}\notag\\
&\le C\sum_{j=0}^{q-1}\norm{e_w^{n-j}}_h
+C(h^{p+1}+\dt^q).
\label{eq:masked_gate_consistency}
\end{align}
The masked gates lie in $[-\varepsilon,1+\varepsilon]$.  The function $G$ and $F$ are
Lipschitz on this fixed compact box.  Applying their pointwise Lipschitz
bounds to \eqref{eq:masked_gate_consistency} proves
\eqref{eq:coefficient_error}.  This argument uses the exact gate bounds,
not bounds on the polynomial interpolants between their nodes.
\end{proof}

Following the standard Ritz projection construction
\cite{thomee2006galerkin}, define $R_h:H^1(\Gamma)\to X_h^p$ by
\begin{equation}\label{eq:ritz_projection}
D(\ds(R_hv-v),\ds\chi_h)_\Gamma+(R_hv-v,\chi_h)_\Gamma=0,
\qquad \chi_h\in X_h^p.
\end{equation}
The second term makes the bilinear form coercive.

\begin{lemma}[Ritz projection on the metric graph]
\label{lem:ritz_graph_projection}
On a family of quasi-uniform edge meshes compatible at the junction,
\begin{equation}\label{eq:ritz_graph_projection_error}
\norm{R_hv-v}+h\norm{\ds(R_hv-v)}
\le Ch^{p+1}\norm{v}_{H^{p+1}(\Gamma)},
\qquad v\in H^{p+1}(\Gamma).
\end{equation}
If $v\in W^{1,\infty}(0,T;H^{p+1}(\Gamma))$, then
\begin{equation}\label{eq:ritz_graph_time_difference}
\frac1{\dt}\norm{\bdf{(R_hv-v)}^{n+1}}
\le Ch^{p+1}
\norm{\partial_tv}_{L^\infty(t_{n+1-q},t_{n+1};H^{p+1}(\Gamma))},
\qquad n\ge q-1.
\end{equation}
\end{lemma}

\begin{proof}
Let
\[
a_\Gamma(u,v):=D(\ds u,\ds v)_\Gamma+(u,v)_\Gamma,
\qquad \rho:=R_hv-v.
\]
Coercivity, Galerkin orthogonality, and
\eqref{eq:graph_interpolation_error} imply
\[
\norm{\rho}_{H^1(\Gamma)}\le Ch^p\norm{v}_{H^{p+1}(\Gamma)}.
\]
For the corresponding metric-graph estimate with linear elements, see
\cite[Theorem~3.2]{arioli2018finite}.
For the $L^2$ estimate, let $\psi\in H^1(\Gamma)$ solve
$a_\Gamma(\psi,\chi)=(\rho,\chi)_\Gamma$ for every $\chi\in H^1(\Gamma)$.
Coercivity gives $\norm{\psi}_{H^1(\Gamma)}\le C\norm{\rho}$.
On each edge, testing with compactly supported functions gives
$-D\ds^2\psi+\psi=\rho$ in the distributional sense.  Hence
$\ds^2\psi\in L^2$ edgewise and
\[
D\norm{\ds^2\psi}\le\norm{\psi}+\norm{\rho}\le C\norm{\rho},
\qquad
\norm{\psi}_{H^2(\Gamma)}\le C\norm{\rho}.
\]
The weak formulation also enforces the natural sealed-end and Kirchhoff
conditions.  In particular, branching does not invalidate this edgewise
regularity estimate.  Since $p\ge1$, interpolation of $\psi$ gives
$\norm{\psi-I_h\psi}_{H^1(\Gamma)}\le Ch\norm{\psi}_{H^2(\Gamma)}$.
Therefore
\[
\norm{\rho}^2=a_\Gamma(\rho,\psi-I_h\psi)
\le Ch\norm{\rho}_{H^1(\Gamma)}\norm{\psi}_{H^2(\Gamma)}
\le Ch\norm{\rho}_{H^1(\Gamma)}\norm{\rho}.
\]
This proves \eqref{eq:ritz_graph_projection_error}.

For the time-difference bound, $R_h$ is linear and independent of time, so
\[
\frac1{\dt}\bdf{(R_hv-v)}^{n+1}
=(R_h-I)\left(\frac1{\dt}\bdf v^{n+1}\right).
\]
The BDF coefficients satisfy $\sum_{j=0}^q a_j^{(q)}=0$.
Consequently,
\[
\frac1{\dt}\bdf v^{n+1}
=\frac1{\dt}\sum_{j=1}^q a_j^{(q)}
\bigl(v(t_{n+1-j})-v(t_{n+1})\bigr),
\]
whose $H^{p+1}(\Gamma)$ norm is at most
\[
\left(\sum_{j=1}^q j|a_j^{(q)}|\right)
\norm{\partial_tv}_{L^\infty(t_{n+1-q},t_{n+1};H^{p+1}(\Gamma))}.
\]
Apply \eqref{eq:ritz_graph_projection_error} to this quotient to obtain
\eqref{eq:ritz_graph_time_difference}.
\end{proof}

Split the voltage error as
\begin{equation}\label{eq:voltage_error_split}
V_h^n-V(t_n)=\theta^n+\rho^n,
\qquad \theta^n:=V_h^n-R_hV(t_n),
\qquad \rho^n:=R_hV(t_n)-V(t_n).
\end{equation}
The preceding lemma and Assumption~\ref{ass:conv_regularity} imply
\begin{equation}\label{eq:ritz_error}
\norm{\rho^n}+h\norm{\ds\rho^n}\le Ch^{p+1},
\qquad
\frac1{\dt}\norm{\bdf\rho^{n+1}}\le Ch^{p+1}.
\end{equation}
The first estimate holds at all time levels; the second holds whenever
the full BDF history is available.

\begin{assumption}[Starting errors]\label{ass:starting_errors}
At the starting levels $j=0,\ldots,q-1$,
\begin{equation}\label{eq:starting_error_bounds}
\norm{V_h^j-V(t_j)}
+\sum_{w\in\mathcal W}\norm{w_h^j-I_hw(t_j)}_h
\le C_{\rm init}(h^{p+1}+\dt^q).
\end{equation}
In addition,
\begin{equation}\label{eq:starting_gradient_bound}
\eta_q\dt^{1/2}\norm{\ds(V_h^{q-1}-R_hV(t_{q-1}))}
\le C_{\rm init}(h^{p+1}+\dt^q).
\end{equation}
\end{assumption}

For $q=1,2$, the chosen multiplier has $\eta_q=0$, so
\eqref{eq:starting_gradient_bound} imposes no additional condition.
For the other orders it controls the gradient contribution to the initial
discrete energy.  It holds, for example, when
$V_h^{q-1}=R_hV(t_{q-1})$.

In the one-step estimates and the main convergence theorem below, all
integrals in the voltage equation are evaluated exactly. 

\subsection{One-step error estimates}

Let $\mathbf M_q$ and $\eta_q$ be the fixed multiplier matrix and parameter in
Lemma~\ref{lem:bdf_multiplier}.  With $\theta^n$ and $e_w^n$ defined in
\eqref{eq:voltage_error_split} and \eqref{eq:gate_error_definition}, respectively,
define the newest-first history vectors for $n\ge q-1$ by
\begin{equation}\label{eq:error_history_vectors}
\bm\Theta^n=(\theta^n,\ldots,\theta^{n-q+1})^T,
\qquad
\bm E_w^n=(e_w^n,\ldots,e_w^{n-q+1})^T.
\end{equation}
The history norms of the vectors in \eqref{eq:error_history_vectors} are those
of Lemma~\ref{lem:bdf_multiplier}, using
$(\cdot,\cdot)_\Gamma$ for $\bm\Theta^n$ and $(\cdot,\cdot)_h$ for
$\bm E_w^n$.  In both cases, \eqref{eq:norm_equiv} applies.
Set
\begin{equation}\label{eq:coupled_error_energies}
\begin{aligned}
A_n&:=\norm{\bm\Theta^n}_{\mathbf M_q}^2
+\frac{D\eta_q^2\dt}{C_m}\norm{\ds\theta^n}^2,
\\
B_n&:=\sum_{w\in\mathcal W}\norm{\bm E_w^n}_{\mathbf M_q,h}^2,
\qquad r^2:=h^{2p+2}+\dt^{2q}.
\end{aligned}
\end{equation}

\begin{lemma}[Voltage error estimate]\label{lem:voltage_error}
Under Assumption~\ref{ass:conv_regularity}, with $A_n$, $B_n$, and $r$
defined in \eqref{eq:coupled_error_energies},
\begin{equation}\label{eq:voltage_error_rec}
A_{n+1}\le(1+C\dt)A_n+C\dt B_n+C\dt r^2,
\qquad n\ge q-1,
\end{equation}
for all steps with $t_{n+1}\le T$.
\end{lemma}

\begin{proof}
The BDF temporal residual of the exact voltage is
\begin{equation}\label{eq:tauV_definition}
\tau_V^{n+1}:=\frac1{\dt}\sum_{j=0}^q a_j^{(q)}V(t_{n+1-j})
-\partial_tV(t_{n+1}).
\end{equation}
By the order-$q$ consistency of the BDF formula, Taylor expansion with an
integral remainder gives
\begin{equation}\label{eq:voltage_time_consistency}
\norm{\tau_V^{n+1}}
\le C\dt^q
\norm{\partial_t^{q+1}V}_{L^\infty(t_{n+1-q},t_{n+1};L^2(\Gamma))}
\le C\dt^q.
\end{equation}
Inserting \eqref{eq:tauV_definition} in the exact weak equation
\eqref{eq:weak_cont} and subtracting it from the exact-integral numerical
equation \eqref{eq:voltage_scheme} yields
\begin{align}
\frac{C_m}{\dt}(\bdf\theta^{n+1},\chi_h)_\Gamma
&+D(\ds\theta^{n+1},\ds\chi_h)_\Gamma
+(G_h^{n+1}\theta^{n+1},\chi_h)_\Gamma
=\mathcal R_V^{n+1}(\chi_h),
\label{eq:voltage_error_equation}
\end{align}
where $G_h^{n+1}$ and $F_h^{n+1}$ are defined in \eqref{def_GF}-\eqref{eq:discrete_coefficients}, $G^{n+1}$ and $F^{n+1}$ are defined in
\eqref{eq:exact_coefficients}, and
\begin{align*}
\mathcal R_V^{n+1}(\chi_h)
={}&-\frac{C_m}{\dt}(\bdf\rho^{n+1},\chi_h)_\Gamma
-D(\ds\rho^{n+1},\ds\chi_h)_\Gamma
-(G_h^{n+1}\rho^{n+1},\chi_h)_\Gamma\\
&-((G_h^{n+1}-G^{n+1})V(t_{n+1}),\chi_h)_\Gamma
+(F_h^{n+1}-F^{n+1},\chi_h)_\Gamma\\
&-C_m(\tau_V^{n+1},\chi_h)_\Gamma.
\end{align*}
Ritz orthogonality gives
\[
D(\ds\rho^{n+1},\ds\chi_h)_\Gamma
=-(\rho^{n+1},\chi_h)_\Gamma.
\]
Using \eqref{eq:ritz_error}, \eqref{eq:voltage_time_consistency},
\eqref{eq:coefficient_error}, the bounded exact voltage, and the uniform
bound $G_h^{n+1}\le g_1$, we obtain
\[
|\mathcal R_V^{n+1}(\chi_h)|
\le C\left(h^{p+1}+\dt^q
+\sum_{w\in\mathcal W}\sum_{j=0}^{q-1}\norm{e_w^{n-j}}_h\right)
\norm{\chi_h}.
\]
The history-norm equivalence in Lemma~\ref{lem:bdf_multiplier} and \eqref{eq:coupled_error_energies} therefore implies
\begin{equation}\label{eq:voltage_residual_bound}
|\mathcal R_V^{n+1}(\chi_h)|
\le C(B_n+r^2)^{1/2}\norm{\chi_h}.
\end{equation}
By Young's inequality,
\[
|\mathcal R_V^{n+1}(\chi_h)|
\le\frac{g_0}2\norm{\chi_h}^2+C(B_n+r^2).
\]
Test \eqref{eq:voltage_error_equation} with $\chi_h=\theta^{n+1}-\eta_q\theta^n$.
The temporal, diffusion, and reaction estimates
\eqref{eq:temporal_est}, \eqref{eq:diff_est}, and \eqref{eq:react_est}
apply with $V_h$ replaced by $\theta$.
After cancellation and multiplication by $2\dt/C_m$, the definitions in
\eqref{eq:coupled_error_energies} give
\begin{align*}
A_{n+1}
&+\frac{D(1-\eta_q^2)\dt}{C_m}\norm{\ds\theta^{n+1}}^2
+\frac{g_0\dt}{C_m}\norm{\theta^{n+1}}^2\\
&\le A_n+\frac{\eta_q^2g_1\dt}{C_m}\norm{\theta^n}^2
+C\dt(B_n+r^2).
\end{align*}
Here a further nonnegative term containing $\norm{\ds\chi_h}^2$ has
been discarded.  Since $\eta_q<1$ and
$\norm{\theta^n}^2\le C A_n$ by
\eqref{eq:norm_equiv}, dropping the remaining
dissipation terms proves \eqref{eq:voltage_error_rec}.
\end{proof}

\begin{lemma}[Gate error estimate]\label{lem:gate_error}
Under Assumption~\ref{ass:conv_regularity}, with $A_n$, $B_n$, and $r$
defined in \eqref{eq:coupled_error_energies},
\begin{equation}\label{eq:gate_error_rec}
B_{n+1}\le(1+C\dt)B_n+C\dt A_{n+1}+C\dt r^2,
\qquad n\ge q-1,
\end{equation}
for all steps with $t_{n+1}\le T$.
\end{lemma}

\begin{proof}
For each gate and each node $a\in\mathcal N_h^p$, define
\begin{equation}\label{eq:gate_tau_definition}
\tau_w^{n+1}(a):=
\frac1{\dt}\sum_{j=0}^q a_j^{(q)}w(a,t_{n+1-j})
-\partial_tw(a,t_{n+1}).
\end{equation}
Taylor expansion and the BDF moment identities give
\[
|\tau_w^{n+1}(a)|
\le C\dt^q
\norm{\partial_t^{q+1}w}_{L^\infty(t_{n+1-q},t_{n+1};C^0(\Gamma))}.
\]
Sobolev embedding, Assumption~\ref{ass:conv_regularity}, and positivity
of the nodal weights therefore yield
\begin{equation}\label{eq:gate_time_consistency}
\norm{\tau_w^{n+1}}_h\le C\dt^q.
\end{equation}
Because the exact voltage is in the uncut range, the exact gate equation
\eqref{eq:gating} uses the same regularized rates.  We subtract nodal equation of \eqref{eq:gating}  from the numerical gate equation \eqref{eq:gate_scheme} to obtain
\begin{equation}\label{eq:gate_error_eq}
\frac1{\dt}\bdf{e_w}^{n+1}
+\widetilde\lambda_w(V_h^{n+1})e_w^{n+1}
=R_w^{n+1},
\end{equation}
as an equality of nodal vectors, where $e_w^{n+1}$ is defined in \eqref{eq:gate_error_definition}, and
\begin{align}
R_w^{n+1}:={}&
\widetilde\alpha_w(V_h^{n+1})
-\widetilde\alpha_w(I_hV(t_{n+1}))\notag\\
&-\bigl[\widetilde\lambda_w(V_h^{n+1})
-\widetilde\lambda_w(I_hV(t_{n+1}))\bigr]w_I^{n+1}
-\tau_w^{n+1}.
\label{eq:gate_residual_definition}
\end{align}
All products and rate evaluations in this identity are nodal.  Since
$I_hV(t_{n+1})(a)=V(a,t_{n+1})$ and
$0\le w_I^{n+1}(a)\le1$, \eqref{eq:rate_lipschitz} implies
\begin{equation}\label{eq:gate_residual_bound}
\norm{R_w^{n+1}}_h
\le C\norm{V_h^{n+1}-I_hV(t_{n+1})}_h+C\dt^q.
\end{equation}
Test \eqref{eq:gate_error_eq} with
$\zeta_w^{n+1}:=e_w^{n+1}-\eta_qe_w^n$ in $(\cdot,\cdot)_h$, and set
$\mathcal E_{e,w}^n:=\norm{\bm E_w^n}_{\mathbf M_q,h}^2$, where
$\bm E_w^n$ is defined in \eqref{eq:error_history_vectors}.
Lemma~\ref{lem:bdf_multiplier}, the weighted square identity and bounds of $\widetilde\lambda_w$ in \eqref{eq:regularized_rate_bounds} give
\begin{align*}
\frac{1}{2\dt}\left(\mathcal E_{e,w}^{n+1}-\mathcal E_{e,w}^n\right)
&+\frac{\lambda_{0}}2\norm{e_w^{n+1}}_h^2
-\frac{\eta_q^2\lambda_{1}}2\norm{e_w^n}_h^2\\
&+\frac{\lambda_{0}}2\norm{\zeta_w^{n+1}}_h^2
\le (R_w^{n+1},\zeta_w^{n+1})_h.
\end{align*}
By Young's inequality,
\[
(R_w^{n+1},\zeta_w^{n+1})_h
\le\frac{1}{2\lambda_{0}}\norm{R_w^{n+1}}_h^2
+\frac{\lambda_{0}}2\norm{\zeta_w^{n+1}}_h^2.
\]
After cancellation, multiplication by $2\dt$, and use of
$c_*\norm{e_w^n}_h^2\le\mathcal E_{e,w}^n$ from
\eqref{eq:norm_equiv} in Lemma~\ref{lem:bdf_multiplier}, we obtain
\[
\mathcal E_{e,w}^{n+1}
\le(1+C\dt)\mathcal E_{e,w}^n
+C\dt\norm{R_w^{n+1}}_h^2.
\]
Summing over the gates and using \eqref{eq:gate_residual_bound} yields
\begin{equation}\label{eq:gate_error_rec_intermediate}
B_{n+1}\le(1+C\dt)B_n
+C\dt\norm{V_h^{n+1}-I_hV(t_{n+1})}_h^2
+C\dt^{2q+1}.
\end{equation}

To estimate the voltage term, write
\[
V_h^{n+1}-I_hV(t_{n+1})
=\theta^{n+1}+R_hV(t_{n+1})-I_hV(t_{n+1}).
\]
Both terms on the right belong to $X_h^p$.
Lemma~\ref{lem:nodal_norm}, the Ritz estimate, and the interpolation
estimate give
\[
\norm{R_hV(t_{n+1})-I_hV(t_{n+1})}_h
\le C\norm{R_hV(t_{n+1})-I_hV(t_{n+1})}
\le Ch^{p+1}.
\]
Since $\norm{\theta^{n+1}}^2\le C A_{n+1}$, it follows that
\begin{equation}\label{eq:voltage_to_gate_bound}
\norm{V_h^{n+1}-I_hV(t_{n+1})}_h^2
\le C A_{n+1}+Ch^{2p+2}.
\end{equation}
Substitution into \eqref{eq:gate_error_rec_intermediate} proves
\eqref{eq:gate_error_rec}.
\end{proof}

\subsection{Coupled convergence}

\begin{theorem}[Convergence under exact integration]
\label{thm:convergence}
Let $q\in\{1,\ldots,5\}$ and an integer $p\ge1$ be fixed, and let the
edge meshes be quasi-uniform and compatible at the junction.  Assume that
all integrals in \eqref{eq:voltage_scheme} are evaluated exactly.  Under
Assumptions~\ref{ass:conv_regularity} and \ref{ass:starting_errors}, there
is a constant $C_{p,q}(T)$ independent of $h,\dt$ such that, for every constant
step $\dt>0$ and every integer
$N\ge q-1$ with $N\dt\le T$,
\begin{equation}\label{eq:convergence_result}
\max_{0\le n\le N}
\left[\norm{V_h^n-V(t_n)}
+\sum_{w\in\mathcal W}\norm{w_h^n-I_hw(t_n)}_h\right]
\le C_{p,q}(T)(h^{p+1}+\dt^q).
\end{equation}
The reconstructed gate functions also satisfy
\begin{equation}\label{eq:continuous_gate_convergence}
\max_{0\le n\le N}\sum_{w\in\mathcal W}
\norm{w_h^n-w(t_n)}
\le C_{p,q}(T)(h^{p+1}+\dt^q).
\end{equation}
No smallness condition on $\dt$ and no relation between $h$ and $\dt$ are
required.
\end{theorem}

\begin{proof}
If $N=q-1$, \eqref{eq:convergence_result} is the starting-error
assumption.  Otherwise, enlarge the constants in
\eqref{eq:voltage_error_rec} and \eqref{eq:gate_error_rec} to one common
constant $C\ge0$.  With $x=C\dt$, the two estimates read
\[
A_{n+1}\le(1+x)A_n+xB_n+xr^2,
\qquad
B_{n+1}\le(1+x)B_n+xA_{n+1}+xr^2.
\]
The order in these recurrences matches the numerical algorithm: the
voltage uses previous gates, and the new gates use the new voltage.
Set $S_n:=A_n+B_n$.  Adding $A_{n+1}$ to both sides of the second inequality, then bounding the right side $A_{n+1}$ by the first inequality,
gives
\begin{align*}
S_{n+1}
&\le(1+x)A_{n+1}+(1+x)B_n+xr^2\\
&\le(1+x)^2A_n+(1+x)^2B_n+(2x+x^2)r^2.
\end{align*}
Equivalently,
\begin{equation}\label{eq:triangular_energy_rec}
S_{n+1}+r^2\le(1+C\dt)^2(S_n+r^2).
\end{equation}
Iterating for $n=q-1,\ldots,m-1$
and using $1+C\dt\le e^{C\dt}$ yields
\begin{equation}\label{eq:coupled_error_energy_bound}
S_m+r^2
\le e^{2C(t_m-t_{q-1})}(S_{q-1}+r^2)
\le e^{2CT}(S_{q-1}+r^2).
\end{equation}
It remains to bound $S_{q-1}$.  The voltage error split
\eqref{eq:voltage_error_split}, the starting-error bound
\eqref{eq:starting_error_bounds}, and the Ritz estimate \eqref{eq:ritz_error}
give, at every starting level,
\[
\norm{\theta^j}
\le\norm{V_h^j-V(t_j)}+\norm{\rho^j}
\le C(h^{p+1}+\dt^q),
\qquad 0\le j\le q-1.
\]
The upper history-norm bound \eqref{eq:norm_equiv} in
Lemma~\ref{lem:bdf_multiplier}, the starting gradient bound
\eqref{eq:starting_gradient_bound}, and the definition of $A_n$ in
\eqref{eq:coupled_error_energies} therefore give
\[
A_{q-1}\le C(h^{2p+2}+\dt^{2q})=Cr^2.
\]
Similarly, the nodal starting-error bound \eqref{eq:starting_error_bounds},
the upper history-norm bound \eqref{eq:norm_equiv}, and the definitions
\eqref{eq:gate_error_definition} and \eqref{eq:coupled_error_energies}
give $B_{q-1}\le Cr^2$.
Thus \eqref{eq:coupled_error_energy_bound} proves
\[
A_m+B_m\le C_{p,q}(T)r^2,
\qquad q-1\le m\le N.
\]
The lower history-norm bound \eqref{eq:norm_equiv}, the definitions
\eqref{eq:gate_error_definition} and \eqref{eq:coupled_error_energies} now imply
\[
\norm{\theta^m}+\sum_{w\in\mathcal W}\norm{e_w^m}_h
\le C_{p,q}(T)(h^{p+1}+\dt^q).
\]
Adding the Ritz projection error bound \eqref{eq:ritz_error} proves
\eqref{eq:convergence_result} at all
computed levels, and \eqref{eq:starting_error_bounds} covers any earlier
starting levels.  Finally, for each gate,
\[
\norm{w_h^n-w(t_n)}
\le C\norm{w_h^n-I_hw(t_n)}_h
+\norm{I_hw(t_n)-w(t_n)}.
\]
The nodal norm equivalence and
\eqref{eq:graph_interpolation_error} prove
\eqref{eq:continuous_gate_convergence}, including the case $N=q-1$.
\end{proof}

\begin{remark}[Extension to finite metric graphs]
\label{rem:general_graph}
The stability and convergence arguments also apply to a fixed connected
metric graph with finitely many edges of positive length.  Impose voltage
continuity and Kirchhoff flux balance at each interior vertex, sealed-end
conditions at degree-one vertices, and the same cable coefficients and gate
kinetics on every edge.  Use compatible quasi-uniform edge meshes and a
single voltage and gate degree of freedom at each vertex.  The assembled
Gauss--Lobatto weights remain positive, the multiplier estimates are
unchanged, and the dual Ritz problem has edgewise $H^2$ regularity with
the standard vertex conditions.  Thus the proofs of
Theorems~\ref{thm:stability}, \ref{thm:gate_stability}, and
\ref{thm:convergence} carry over under the corresponding regularity and
starting-value assumptions.  Their constants may depend on the fixed
graph, including its edge lengths and connectivity.
\end{remark}

\subsection{Extension to numerical quadrature}

Let the voltage equation now be evaluated by a numerical quadrature rule.
At the computed solution, define its full exact-form defect by
\begin{align}
\mathcal Q_h^{n+1}(\chi_h)
:={}&\frac{C_m}{\dt}(\bdf{V_h}^{n+1},\chi_h)_\Gamma
+D(\ds V_h^{n+1},\ds\chi_h)_\Gamma\notag\\
&+(G_h^{n+1}V_h^{n+1}-F_h^{n+1}-I_{\rm ext}^{n+1},\chi_h)_\Gamma,
\qquad \chi_h\in X_h^p,
\label{eq:quadrature_defect_definition}
\end{align}
where every integral on the right is exact.  Thus
$\mathcal Q_h^{n+1}$ is the residual obtained when the quadrature-based
solution is inserted into the exact-integral voltage equation.

\begin{assumption}[Quadrature defect bound]\label{ass:quadrature}
There is a constant $C_Q$, independent of $h$ and $\dt$, such that
\begin{equation}\label{eq:quadrature_consistency}
\mathfrak q_{n+1}:=
\sup_{\chi_h\in X_h^p,\ \chi_h\ne0}
\frac{|\mathcal Q_h^{n+1}(\chi_h)|}{\norm{\chi_h}}
\le C_Qh^{p+1},
\qquad n\ge q-1.
\end{equation}
\end{assumption}

\begin{corollary}[Convergence with numerical quadrature]
\label{cor:quadrature_convergence}
Let the edge meshes be quasi-uniform and compatible at the junction,
and use a constant time step as in Theorem~\ref{thm:convergence}.
Under Assumptions~\ref{ass:conv_regularity}, \ref{ass:starting_errors},
and \ref{ass:quadrature}, the conclusions
\eqref{eq:convergence_result}--\eqref{eq:continuous_gate_convergence}
hold for the quadrature-based scheme, with $C_{p,q}(T)$ also allowed to depend
on $C_Q$.
\end{corollary}

\begin{proof}
The voltage error equation \eqref{eq:voltage_error_equation} gains the
additional term $\mathcal Q_h^{n+1}(\chi_h)$ on its right-hand side.
Consequently,
\[
|\mathcal R_V^{n+1}(\chi_h)+\mathcal Q_h^{n+1}(\chi_h)|
\le\left[C(B_n+r^2)^{1/2}+\mathfrak q_{n+1}\right]
\norm{\chi_h}.
\]
The same multiplier argument and Young's inequality therefore replace
\eqref{eq:voltage_error_rec} by
\[
A_{n+1}\le(1+C\dt)A_n+C\dt B_n
+C\dt\left(r^2+\mathfrak q_{n+1}^2\right).
\]
By \eqref{eq:quadrature_consistency},
$\mathfrak q_{n+1}^2\le C_Q^2h^{2p+2}\le C r^2$.
Thus the original voltage recurrence is recovered with a modified
constant.  The gate recurrence is unchanged, so the proof of
Theorem~\ref{thm:convergence} applies verbatim.
\end{proof}

\begin{remark}
Exact integration gives $\mathcal Q_h^{n+1}=0$.
For a numerical quadrature rule, verifying
\eqref{eq:quadrature_consistency} requires a separate analysis.
\end{remark}

\section{Numerical Experiments}\label{sec:numerical_experiments}

	Unless stated otherwise, the experiments use the Y-shaped graph in
	Figure~\ref{fig:axon_geometry}, with three unit-length branches.
	Section~\ref{subsec:multiple_junctions} also considers a tree with two
	junctions and unequal edge lengths.  All experiments use the physical
	parameters in Table~\ref{tab:parameters}.  The computations use $p=1,\ldots,4$.
	The gate mask and voltage cutoff use
	$\varepsilon=10^{-3}$ and $M=200$, respectively.  The mask parameter
	satisfies $g_{\mathrm{Na}}\varepsilon(1+\varepsilon)^3\le g_{\mathrm L}/2$.
	The forcing, starting histories, final times, and refinement choices are
	specified separately for each experiment below.
	At a grid time $t=t_n$, we write $V_h(t)=V_h^n$ and $w_h(t)=w_h^n$.
	All sampled times below lie on the corresponding time grids; the same
	convention applies to the numerical reference solutions on their own grids.

	On each edge $A_i$ of length $L_i$, let $\ell_i\in[0,L_i]$ denote
	arc length measured from the initial vertex of the oriented edge.
	For the unit-length Y-shaped graph, $L_i=1$ and
	\begin{equation}\label{l-i-s}
	\ell_0(s)=s+1,\qquad \ell_1(s)=\ell_2(s)=s.
	\end{equation}

	\begin{table}[htbp]
		\centering
		\caption{Physical parameters shared by the numerical experiments.}
		\label{tab:parameters}
		\begin{tabular}{@{}llll@{}}
			\toprule
			\textbf{Parameter} & \textbf{Value} & \textbf{Unit} & \textbf{Description} \\
			\midrule
			$C_m$ & $1.0$ & $\mu\mathrm{F}/\mathrm{cm}^2$ & Specific membrane capacitance \\
			$a$ & $0.1$ & $\mu\mathrm{m}$ & Axon radius \\
			$R$ & $35.4$ & $\Omega\,\mathrm{cm}$ & Intracellular resistivity \\
			$g_{\mathrm{Na}}$ & $120.0$ & $\mathrm{mS}/\mathrm{cm}^2$ & Maximum Na$^+$ conductance \\
			$g_{\mathrm{K}}$ & $36.0$ & $\mathrm{mS}/\mathrm{cm}^2$ & Maximum K$^+$ conductance \\
			$g_{\mathrm{L}}$ & $0.3$ & $\mathrm{mS}/\mathrm{cm}^2$ & Maximum leak conductance \\
			$E_{\mathrm{Na}}$ & $115.0$ & $\mathrm{mV}$ & Na$^+$ reversal potential \\
			$E_{\mathrm{K}}$ & $-12.0$ & $\mathrm{mV}$ & K$^+$ reversal potential \\
			$E_{\mathrm{L}}$ & $10.5989\ldots $ & $\mathrm{mV}$ & Leakage reversal potential \\
			\bottomrule
		\end{tabular}
	\end{table}
\noindent We choose $a=0.1\,\mu\mathrm{m}$, consistent in scale with thin
parallel fibers~\cite{sultan2000diameter}.  The electrical parameters follow
Hodgkin and Huxley~\cite{hodgkin1952quantitative}, with $E_{\mathrm L}$
adjusted to enforce the exact resting equilibrium.
Define the equilibrium gate functions by
	\[
	w_\infty(V):=\frac{\alpha_w(V)}{\alpha_w(V)+\beta_w(V)},
	\qquad w\in\mathcal W.
	\]
	In all experiments, the leak reversal potential is chosen as
	\begin{equation}\label{eq:resting_leak_reversal}
		E_{\mathrm{L}}^*
		=
		\frac{
			g_{\mathrm{Na}}m_\infty(0)^3h_\infty(0)(0-E_{\mathrm{Na}})
			+
			g_{\mathrm{K}}n_\infty(0)^4(0-E_{\mathrm{K}})
		}{g_{\mathrm{L}}}
		=
		10.598920969\ldots
	\end{equation}
	The computations use the double-precision value $E_{\mathrm{L}}=E_{\mathrm{L}}^*$;
	Table~\ref{tab:parameters} displays only a rounded decimal.  This choice
	makes $V=0$ with $w=w_\infty(0)$ an exact equilibrium of the unforced
	model and supplies the resting histories in
	Section~\ref{subsec:computational_setup}.

	\subsection{Hodgkin--Huxley Model with a Manufactured Solution}
	\label{subsec:manufactured_solution}
\begin{figure}[htbp]
		\centering
		\includegraphics[width=\textwidth]{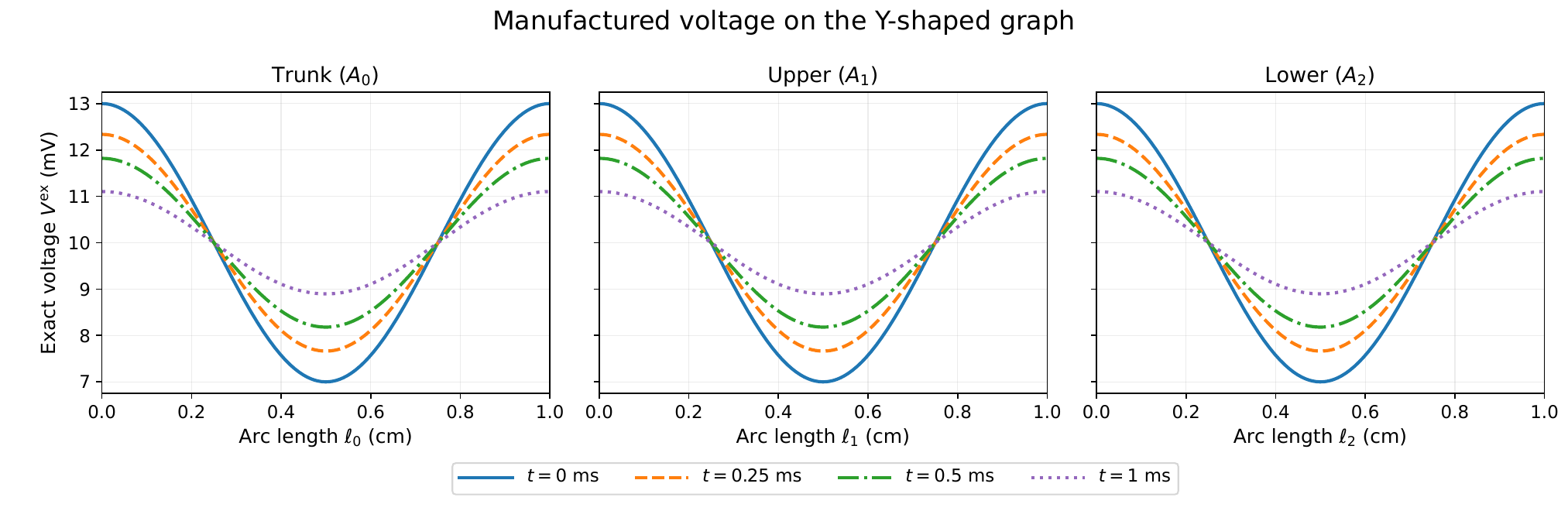}
		\caption{Exact manufactured voltage
		on the trunk and the two daughter branches.}
		\label{fig:manufactured_voltage_profiles}
	\end{figure}

	We first verify the numerical method using prescribed smooth voltage and
	gate fields, with source terms chosen so that these fields solve the
	coupled equations exactly.  We consider the four pairs
	$(q,p)=(2,1),(3,2),(4,3),(5,4)$ in joint space--time refinement.
	Using \eqref{l-i-s}, define
	\begin{equation}\label{eq:manufactured_profile}
		Q_i(s,t):=e^{-t}\cos\bigl(\kappa_i\ell_i(s)\bigr),
		\qquad \kappa_i:=\frac{2\pi}{L_i}.
	\end{equation}
	The manufactured solution is
	\begin{equation}\label{eq:manufactured_fields}
		\begin{aligned}
			V_i^{\mathrm{ex}}(s,t)&=10+3Q_i(s,t),\\
			m_i^{\mathrm{ex}}(s,t)&=0.1+0.02Q_i(s,t),\\
			h_i^{\mathrm{ex}}(s,t)&=0.6+0.03Q_i(s,t),\\
			n_i^{\mathrm{ex}}(s,t)&=0.3+0.02Q_i(s,t).
		\end{aligned}
	\end{equation}
	For $t\ge0$, the gate mask and the voltage cutoff act as the identity
	on the exact fields. At every edge endpoint, $Q_i=e^{-t}$ and $\ds Q_i=0$.
	Thus the fields have a common junction trace, and the voltage satisfies
	the sealed-end and Kirchhoff conditions.
    Figure~\ref{fig:manufactured_voltage_profiles} shows the same cosine
	voltage profile on each branch.

	For each order $q$, the starting levels use the manufactured histories:
	\[
	V_h^j=R_hV^{\mathrm{ex}}(t_j),\qquad
	w_h^j=I_hw^{\mathrm{ex}}(t_j),
	\qquad j=0,\ldots,q-1,\quad w\in\mathcal W,
	\]
	with the Ritz projection defined in \eqref{eq:ritz_projection}.
	Nonpolynomial integrals are evaluated using $24$-point Gauss--Legendre
	quadrature on each element, and the continuous error norms use a
	$40$-point rule.

    \paragraph{Joint space--time refinement.}
    We set $T=1$ and refine space and time together according to
    \[
    h=\dt=\frac1{N_e},\qquad N_e=40,80,160,320,
    \]
    where $N_e$ is the number of elements on each unit-length branch.
    Since $q=p+1$, the bound $h^{p+1}+\dt^q$ predicts joint order $q$.
    For these tests, both the voltage and the reconstructed gates are
    compared directly with the analytic continuous fields:
    \begin{equation}\label{eq:manufactured_joint_errors}
        \begin{aligned}
            E_V^{\mathrm{joint}}
            &:=\max_{t\in\mathcal T}
              \norm{V_h(t)-V^{\mathrm{ex}}(t)},\\
            E_W^{\mathrm{joint}}
            &:=\max_{t\in\mathcal T}\max_{w\in\mathcal W}
              \norm{w_h(t)-w^{\mathrm{ex}}(t)},\\
            \mathcal T&:=\{0,0.05,\ldots,1\}.
        \end{aligned}
    \end{equation}
    In particular, the joint gate error includes the continuous
    interpolation error. Rates are computed as
    $\log_2(E_{\mathrm{coarse}}/E_{\mathrm{fine}})$.
    Tables~\ref{tab:manufactured_joint_bdf2_p1}--\ref{tab:manufactured_joint_bdf4_p3}
    report double-precision results with rates approaching the predicted orders.
    For BDF5--$P_4$, we use the \texttt{mpmath} library~\cite{mpmath2023} with $50$ decimal
    digits to prevent floating-point roundoff from obscuring fifth-order
    convergence, see Table~\ref{tab:manufactured_joint_bdf5_p4}.

    \begin{table}[htbp]
        \centering
        \caption{Joint refinement for BDF2--$P_1$ with $T=1$ and $h=\dt$.}
        \label{tab:manufactured_joint_bdf2_p1}
        \begin{tabular}{@{}r r r r r r@{}}
            \toprule
            $N_e$ & $h=\dt$ & $E_V^{\mathrm{joint}}$ & Rate
            & $E_W^{\mathrm{joint}}$ & Rate \\
            \midrule
            $40$ & $1/40$ & $3.3888\times10^{-3}$ & ---
            & $8.2715\times10^{-5}$ & --- \\
            $80$ & $1/80$ & $8.4534\times10^{-4}$ & $2.0032$
            & $2.0687\times10^{-5}$ & $1.9994$ \\
            $160$ & $1/160$ & $2.1122\times10^{-4}$ & $2.0008$
            & $5.1723\times10^{-6}$ & $1.9999$ \\
            $320$ & $1/320$ & $5.2798\times10^{-5}$ & $2.0002$
            & $1.2931\times10^{-6}$ & $2.0000$ \\
            \bottomrule
        \end{tabular}
    \end{table}

    \begin{table}[htbp]
        \centering
        \caption{Joint refinement for BDF3--$P_2$ with $T=1$ and $h=\dt$.}
        \label{tab:manufactured_joint_bdf3_p2}
        \begin{tabular}{@{}r r r r r r@{}}
            \toprule
            $N_e$ & $h=\dt$ & $E_V^{\mathrm{joint}}$ & Rate
            & $E_W^{\mathrm{joint}}$ & Rate \\
            \midrule
            $40$ & $1/40$ & $8.4024\times10^{-5}$ & ---
            & $8.1862\times10^{-7}$ & --- \\
            $80$ & $1/80$ & $1.1092\times10^{-5}$ & $2.9213$
            & $1.0235\times10^{-7}$ & $2.9996$ \\
            $160$ & $1/160$ & $1.4249\times10^{-6}$ & $2.9605$
            & $1.2795\times10^{-8}$ & $2.9999$ \\
            $320$ & $1/320$ & $1.8057\times10^{-7}$ & $2.9803$
            & $1.5994\times10^{-9}$ & $3.0000$ \\
            \bottomrule
        \end{tabular}
    \end{table}

    \begin{table}[htbp]
        \centering
        \caption{Joint refinement for BDF4--$P_3$ with $T=1$ and $h=\dt$.}
        \label{tab:manufactured_joint_bdf4_p3}
        \begin{tabular}{@{}r r r r r r@{}}
            \toprule
            $N_e$ & $h=\dt$ & $E_V^{\mathrm{joint}}$ & Rate
            & $E_W^{\mathrm{joint}}$ & Rate \\
            \midrule
            $40$ & $1/40$ & $1.7762\times10^{-6}$ & ---
            & $1.4915\times10^{-8}$ & --- \\
            $80$ & $1/80$ & $1.2208\times10^{-7}$ & $3.8629$
            & $1.0371\times10^{-9}$ & $3.8462$ \\
            $160$ & $1/160$ & $8.0091\times10^{-9}$ & $3.9300$
            & $6.8419\times10^{-11}$ & $3.9220$ \\
            $320$ & $1/320$ & $5.1763\times10^{-10}$ & $3.9516$
            & $4.4372\times10^{-12}$ & $3.9467$ \\
            \bottomrule
        \end{tabular}
    \end{table}

    \begin{table}[htbp]
        \centering
        \caption{Joint refinement for BDF5--$P_4$ with $T=1$ and $h=\dt$.}
        \label{tab:manufactured_joint_bdf5_p4}
        \begin{tabular}{@{}r r r r r r@{}}
            \toprule
            $N_e$ & $h=\dt$ & $E_V^{\mathrm{joint}}$ & Rate
            & $E_W^{\mathrm{joint}}$ & Rate \\
            \midrule
            $40$ & $1/40$ & $4.0277\times10^{-8}$ & ---
            & $3.3569\times10^{-10}$ & --- \\
            $80$ & $1/80$ & $1.4300\times10^{-9}$ & $4.8158$
            & $1.2118\times10^{-11}$ & $4.7920$ \\
            $160$ & $1/160$ & $4.7645\times10^{-11}$ & $4.9076$
            & $4.0671\times10^{-13}$ & $4.8970$ \\
            $320$ & $1/320$ & $1.5375\times10^{-12}$ & $4.9537$
            & $1.3170\times10^{-14}$ & $4.9487$ \\
            \bottomrule
        \end{tabular}
    \end{table}

	\subsection{Multiple Junctions and Unequal Edge Lengths}
	\label{subsec:multiple_junctions}
\begin{figure}[htbp]
		\centering
		\includegraphics[width=0.88\textwidth]{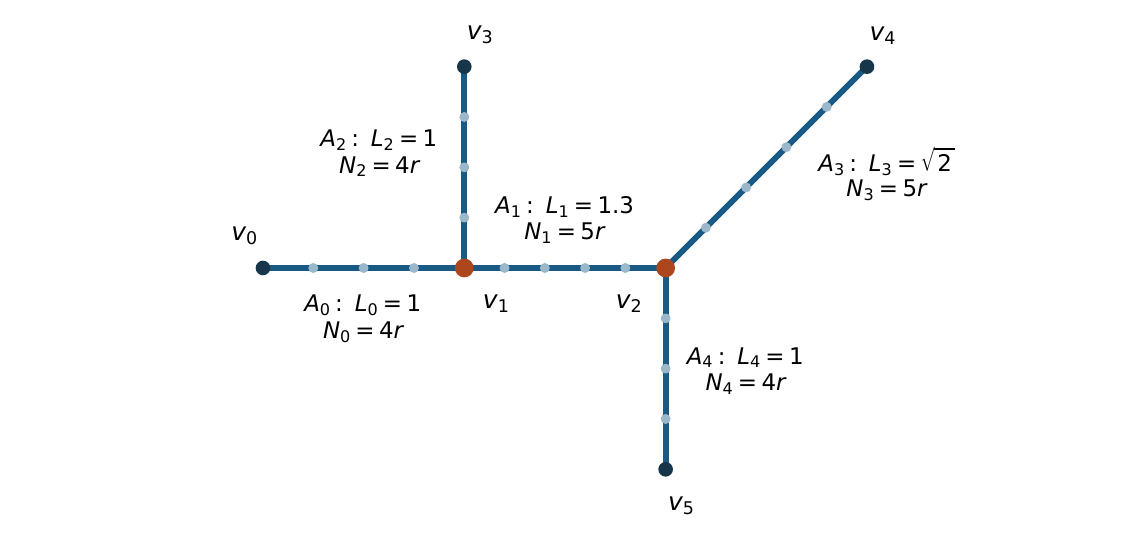}
		\caption{Tree with two junctions and unequal edge lengths.
		Orange vertices are junctions.  Labels give the metric lengths and
		element counts; small light dots mark element boundaries for $r=1$.}
		\label{fig:two_junction_geometry}
	\end{figure}

	To check the assembly and spatial approximation beyond a single
	bifurcation, we use the five-edge tree in
	Figure~\ref{fig:two_junction_geometry}.  Its oriented edges are
	\[
	A_0=(v_0,v_1),\quad A_1=(v_1,v_2),\quad A_2=(v_1,v_3),\quad
	A_3=(v_2,v_4),\quad A_4=(v_2,v_5),
	\]
	with metric lengths $(L_0,\ldots,L_4)=(1,1.3,1,\sqrt{2},1)$.
	Both $v_1$ and $v_2$ are junctions of degree three.  A single global
	degree of freedom is shared by all incident edges at each junction;
	sealed conditions are imposed at the four leaves, and Kirchhoff balance
	is incorporated through the assembled weak formulation.

\begin{figure}[htbp]
		\centering
		\includegraphics[width=\textwidth]{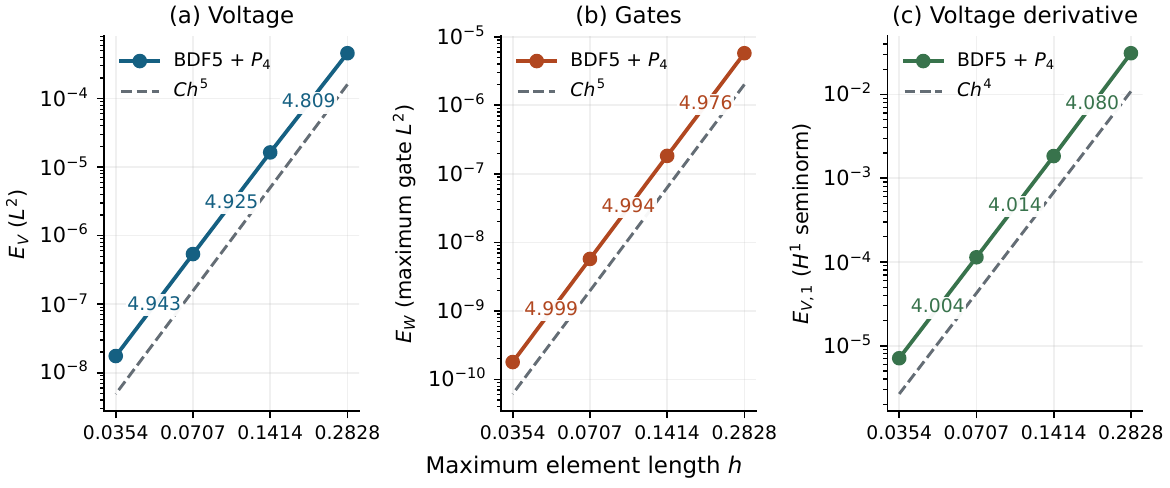}
		\caption{Spatial convergence on the two-junction tree with BDF5--$P_4$
		and fixed $\dt=0.005$: (a) voltage $L^2$ error, (b) maximum gate $L^2$
		error, and (c) voltage derivative error ($H^1$ seminorm).  Errors are
		maxima over $\mathcal T_g$.  Dashed lines indicate fifth-order scaling
		in (a)--(b) and fourth-order scaling in (c); numbers give observed rates
		between successive meshes.}
		\label{fig:two_junction_convergence}
	\end{figure}

	We reuse the manufactured fields in
	\eqref{eq:manufactured_profile}-\eqref{eq:manufactured_fields}, with
	$\kappa_i=2\pi/L_i$ on each edge.  At every vertex these fields have
	the same trace and zero edgewise derivative, so they satisfy the vertex
	and leaf conditions.  BDF5--$P_4$ is used with $T=0.1$, fixed
	$\dt=0.005$, and the manufactured starting histories specified in
	Section~\ref{subsec:manufactured_solution}.  The element counts are
	\[
	(N_0,\ldots,N_4)=r(4,5,4,5,4),\qquad r=1,2,4,8,
	\]
	so $h=\sqrt{2}/(5r)$ and the global number of degrees of freedom is
	$88r+1$.  Assembly uses $24$-point Gauss--Legendre quadrature and
	continuous error evaluation uses a $40$-point rule.  At the sampled
	times $\mathcal T_g=\{0,0.025,0.05,0.075,0.1\}$, define
	\begin{equation}\label{eq:multiple_junction_errors}
	\begin{aligned}
	E_V&:=\max_{t\in\mathcal T_g}\norm{V_h(t)-V^{\rm ex}(t)},\\
	E_W&:=\max_{t\in\mathcal T_g}\max_{w\in\mathcal W}
	       \norm{w_h(t)-w^{\rm ex}(t)},\\
	E_{V,1}&:=\max_{t\in\mathcal T_g}
	       \norm{\ds(V_h(t)-V^{\rm ex}(t))}.
	\end{aligned}
	\end{equation}
	Here $E_W$ includes gate interpolation error, and $E_{V,1}$ is the
	voltage $H^1$ seminorm error.
	Figure~\ref{fig:two_junction_convergence} shows $L^2$ rates approaching
	five for voltage and gates and an $H^1$ seminorm rate approaching four
	for voltage, as expected for $P_4$.

	\subsection{Hodgkin--Huxley Model with a Localized Stimulus}\label{subsec:computational_setup}
	The external stimulus acts only on the trunk
	$A_0$, with amplitude $I_0=100.0$:
	\begingroup
	\allowdisplaybreaks[1]
	\begin{align}
			I_{\mathrm{ext},0}(s,t)
			&=
			I_0\,\phi_t(t)\phi_s(s),
			\qquad s\in[-1,0],\notag
			\\
			I_{\mathrm{ext},1}(s,t)
			&=
			I_{\mathrm{ext},2}(s,t)=0,\notag
			\\[0.5em]
			\phi_t(t)
			&=
			\begin{cases}
				e\exp\!\left(-\dfrac{0.01}{0.01-(t-0.2)^2}\right),
				& t\in(0.1,0.3),\\[0.5em]
				0, & \text{otherwise},
			\end{cases}\notag
			\\[0.5em]
			\phi_s(s)
			&=
			\begin{cases}
				e\exp\!\left(-\dfrac{0.01}{0.01-(2s+1.9)^2}\right),
				& s\in(-1.0,-0.9),\\[0.5em]
				0, & \text{otherwise}.
			\end{cases}
	\label{eq:stimulus_definition}
	\end{align}
	\endgroup
	This stimulus is smooth and compactly supported in time and space, so it
	satisfies \eqref{eq:external_stimulus_bound}.

\begin{figure}[htp]
		\centering
		\includegraphics[width=\textwidth]{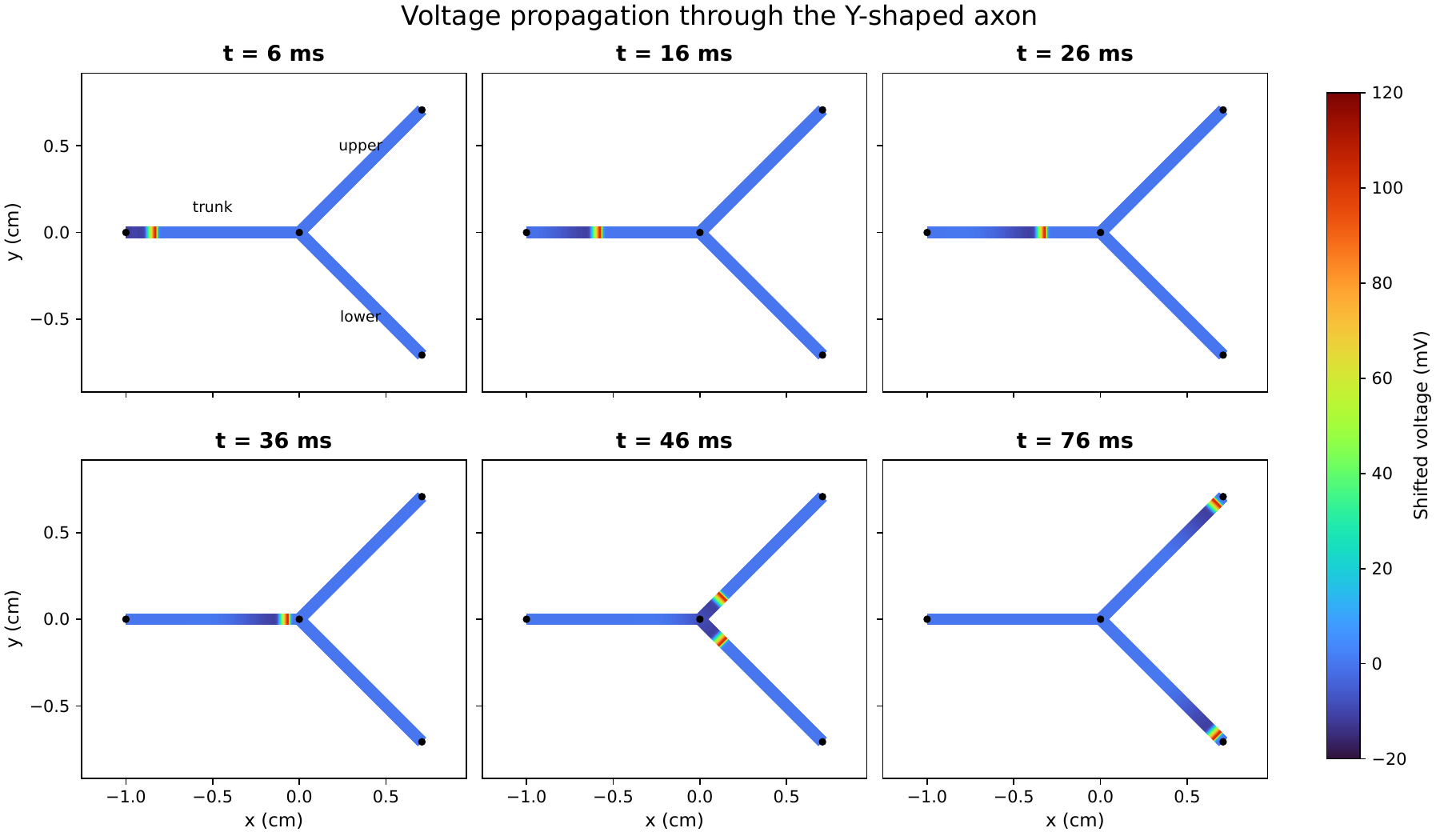}
		\caption{Voltage propagation under the localized stimulus, computed
		with BDF3--$P_2$, $h=0.01\,\mathrm{cm}$, and
		$\dt=0.02\,\mathrm{ms}$.  The panels show
		$t=6,16,26,36,46,$ and $76\,\mathrm{ms}$, from left to right and
		top to bottom.  The shared color scale gives the shifted voltage in
		$\mathrm{mV}$.}
		\label{fig:localized_voltage_propagation}
	\end{figure}

	We initialize at the resting equilibrium $V^0=0$, $w^0=w_\infty(0)$.
	Since the stimulus starts at $t=0.1\,\mathrm{ms}$ and all simulations
	satisfy $(q-1)\dt<0.1\,\mathrm{ms}$, every starting level is an exact
	resting state:
	\begin{equation}\label{eq:initialization}
		V_h^j=0,
		\qquad
		w_h^j=w_\infty(0),
		\qquad j=0,\ldots,q-1,\quad w\in\mathcal{W}.
	\end{equation}
	These constants have zero starting and Ritz projection errors and
	satisfy Assumption~\ref{ass:starting_errors}.
	The corresponding resting gate values are
	\[
		m_\infty(0)\approx0.05293,
		\qquad
		h_\infty(0)\approx0.59612,
		\qquad
		n_\infty(0)\approx0.31768.
	\]
	\paragraph{Voltage propagation.}
	We use a final time of $T=76\,\mathrm{ms}$.
	Figure~\ref{fig:localized_voltage_propagation} shows the pulse traveling
	along the trunk and splitting symmetrically into the two daughter branches.

\begin{figure}[htbp]
		\centering
		\includegraphics[width=\textwidth]{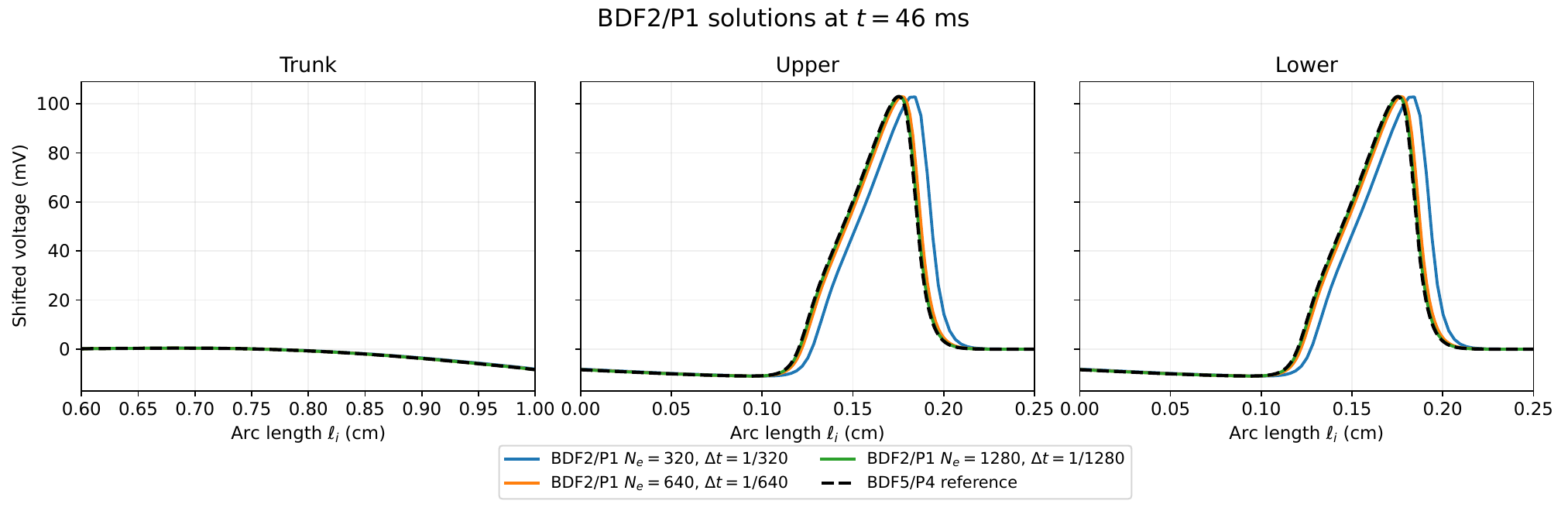}
		\caption{Zoomed BDF2--$P_1$ voltage profiles at $t=46\,\mathrm{ms}$
		for $N_e=320,640,1280$, with $h=1/N_e\,\mathrm{cm}$ and
		$\dt=1/N_e\,\mathrm{ms}$.  The dashed curve is the common
		BDF5--$P_4$ reference used in Table~\ref{tab:physical_joint_convergence}.}
		\label{fig:bdf2_p1_profiles_t46}
	\end{figure}

	\paragraph{Joint space--time convergence.}
	We also compare the physical Hodgkin--Huxley fields at
	$t_*=46\,\mathrm{ms}$ under joint refinement, with
	$h=1/N_e\,\mathrm{cm}$ and $\dt=1/N_e\,\mathrm{ms}$.
	All comparisons use the same BDF5--$P_4$ numerical reference with
	$1280$ elements per branch and $\dt_{\rm ref}=1/1280\,\mathrm{ms}$.
	Assembly uses $24$-point Gauss--Legendre quadrature.  Define the
	continuous graph errors
	\begin{equation}\label{eq:physical_joint_errors}
		\begin{aligned}
		E_V^{46}&:=\norm{V_h(t_*)-V_{\rm ref}(t_*)},\\
		E_W^{46}&:=\max_{w\in\mathcal W}
		\norm{w_h(t_*)-w_{\rm ref}(t_*)},\\
		E_{V,H^1}^{46}&:=\left[(E_V^{46})^2+
		\norm{\ds\bigl(V_h(t_*)-V_{\rm ref}(t_*)\bigr)}^2\right]^{1/2}.
		\end{aligned}
	\end{equation}
	Table~\ref{tab:physical_joint_convergence} reports these errors and the
	rates $\log_2(E_{N_e}/E_{2N_e})$.  The finest BDF2--$P_1$ $L^2$ rates are
	close to two, while the BDF4--$P_3$ and BDF5--$P_4$ rates approach four
	and five.  The BDF3--$P_2$ $L^2$ slopes remain above three on the available
	grids.  The finest voltage $H^1$ rates for $P_2$, $P_3$, and $P_4$
	are close to two, three, and four, respectively; the $P_1$ rate has not
	yet settled to one.  These are comparisons with a finite numerical
	reference; above-nominal slopes do not establish a higher formal order.
	Figure~\ref{fig:bdf2_p1_profiles_t46} illustrates the approach of the
	BDF2--$P_1$ pulse profiles to the reference under refinement.

	\begin{table}[htbp]
		\centering
		\caption{Joint refinement for the physical Hodgkin--Huxley problem at
		$t=46\,\mathrm{ms}$.}
		\label{tab:physical_joint_convergence}
		\begin{tabular}{@{}r r r r r r r@{}}
			\toprule
			$N_e$ & $E_V^{46}$ & Rate & $E_W^{46}$ & Rate & $E_{V,H^1}^{46}$ & Rate \\
			\midrule
			\multicolumn{7}{c}{BDF2--$P_1$} \\
			\midrule
			$80$ & $3.9162\times10^{1}$ & --- & $4.4531\times10^{-1}$ & --- & $1.7855\times10^{3}$ & --- \\
			$160$ & $2.8259\times10^{1}$ & $0.4707$ & $2.9816\times10^{-1}$ & $0.5787$ & $1.8869\times10^{3}$ & $-0.0797$ \\
			$320$ & $9.3142\times10^{0}$ & $1.6012$ & $9.3383\times10^{-2}$ & $1.6749$ & $1.1639\times10^{3}$ & $0.6971$ \\
			$640$ & $2.3901\times10^{0}$ & $1.9623$ & $2.3940\times10^{-2}$ & $1.9637$ & $3.4881\times10^{2}$ & $1.7384$ \\
			$1280$ & $5.9685\times10^{-1}$ & $2.0017$ & $5.9793\times10^{-3}$ & $2.0014$ & $9.5080\times10^{1}$ & $1.8752$ \\
			\midrule
			\multicolumn{7}{c}{BDF3--$P_2$} \\
			\midrule
			$80$ & $2.1457\times10^{1}$ & --- & $2.2041\times10^{-1}$ & --- & $1.7926\times10^{3}$ & --- \\
			$160$ & $2.3469\times10^{0}$ & $3.1926$ & $2.3495\times10^{-2}$ & $3.2298$ & $3.4564\times10^{2}$ & $2.3747$ \\
			$320$ & $1.6159\times10^{-1}$ & $3.8603$ & $1.6167\times10^{-3}$ & $3.8612$ & $3.0136\times10^{1}$ & $3.5197$ \\
			$640$ & $1.0437\times10^{-2}$ & $3.9526$ & $1.0447\times10^{-4}$ & $3.9519$ & $5.0169\times10^{0}$ & $2.5866$ \\
			$1280$ & $6.7427\times10^{-4}$ & $3.9522$ & $6.7640\times10^{-6}$ & $3.9490$ & $1.2063\times10^{0}$ & $2.0562$ \\
			$2560$ & $4.5786\times10^{-5}$ & $3.8803$ & $4.6276\times10^{-7}$ & $3.8695$ & $3.0174\times10^{-1}$ & $1.9993$ \\
			\midrule
			\multicolumn{7}{c}{BDF4--$P_3$} \\
			\midrule
			$80$ & $4.5711\times10^{0}$ & --- & $4.5705\times10^{-2}$ & --- & $6.4449\times10^{2}$ & --- \\
			$160$ & $1.2888\times10^{-1}$ & $5.1484$ & $1.2884\times10^{-3}$ & $5.1487$ & $2.1321\times10^{1}$ & $4.9178$ \\
			$320$ & $2.4644\times10^{-3}$ & $5.7087$ & $2.4902\times10^{-5}$ & $5.6932$ & $2.0798\times10^{0}$ & $3.3578$ \\
			$640$ & $6.2218\times10^{-5}$ & $5.3077$ & $6.5536\times10^{-7}$ & $5.2478$ & $2.8153\times10^{-1}$ & $2.8851$ \\
			$1280$ & $3.0331\times10^{-6}$ & $4.3585$ & $3.3021\times10^{-8}$ & $4.3108$ & $3.5440\times10^{-2}$ & $2.9898$ \\
			\midrule
			\multicolumn{7}{c}{BDF5--$P_4$} \\
			\midrule
			$80$ & $7.1823\times10^{-1}$ & --- & $7.1537\times10^{-3}$ & --- & $1.0780\times10^{2}$ & --- \\
			$160$ & $6.7239\times10^{-3}$ & $6.7390$ & $6.7835\times10^{-5}$ & $6.7205$ & $4.0406\times10^{0}$ & $4.7376$ \\
			$320$ & $7.6875\times10^{-5}$ & $6.4507$ & $7.8463\times10^{-7}$ & $6.4339$ & $2.9298\times10^{-1}$ & $3.7857$ \\
			$640$ & $2.0972\times10^{-6}$ & $5.1960$ & $2.3119\times10^{-8}$ & $5.0849$ & $1.6429\times10^{-2}$ & $4.1565$ \\
			\bottomrule
		\end{tabular}
	\end{table}

	\subsection{Runtime Benchmark}\label{subsec:runtime_benchmark}
	We compare computational work at a prescribed solution error, following
	the accuracy--cost perspective used in electrophysiological
	time-integration studies \cite{roy2020timestepping,lindner2023splitting}.
	Here the comparison is between members of the same BDF--FEM family
	under common implementation and reference-solution settings.
	At $T=46\,\mathrm{ms}$, we compare accuracy and computational cost for the
	same problem as in Section~\ref{subsec:computational_setup}, using the same physical
	parameters and resting histories.
	The pairs BDF2--$P_1$, BDF3--$P_2$, and BDF4--$P_3$ use
	$N_e=80,160,320,640,1280$ elements per branch; BDF5--$P_4$ uses
	$N_e=80,160,320,640$.  In every case,
	$h=1/N_e\,\mathrm{cm}$ and $\dt=1/N_e\,\mathrm{ms}$.
	All cases use an optimized evaluation of the gate mask $\phi$ from
	Appendix~\ref{app:masking_function}.  If all extrapolated gate values at
	the assembly quadrature points lie in $[0,1]$, we return a copy of these
	values directly, since $\phi(x)=x$ there.  This avoids evaluating the
	polynomial and hyperbolic-tangent saturation branches.  If any value
	lies outside $[0,1]$, the full mask is evaluated on the entire array.
	This shortcut preserves the mask values and the numerical scheme.
	Assembly uses $24$-point Gauss--Legendre quadrature.  The errors in
	\eqref{eq:physical_joint_errors} are evaluated against the common
	BDF5--$P_4$ reference with $N_e=1280$ and
	$\dt_{\rm ref}=1/1280\,\mathrm{ms}$, using $8$-point quadrature on the
	union of the candidate and reference element boundaries.

\begin{figure}[htbp]
		\centering
		\includegraphics[width=\textwidth]{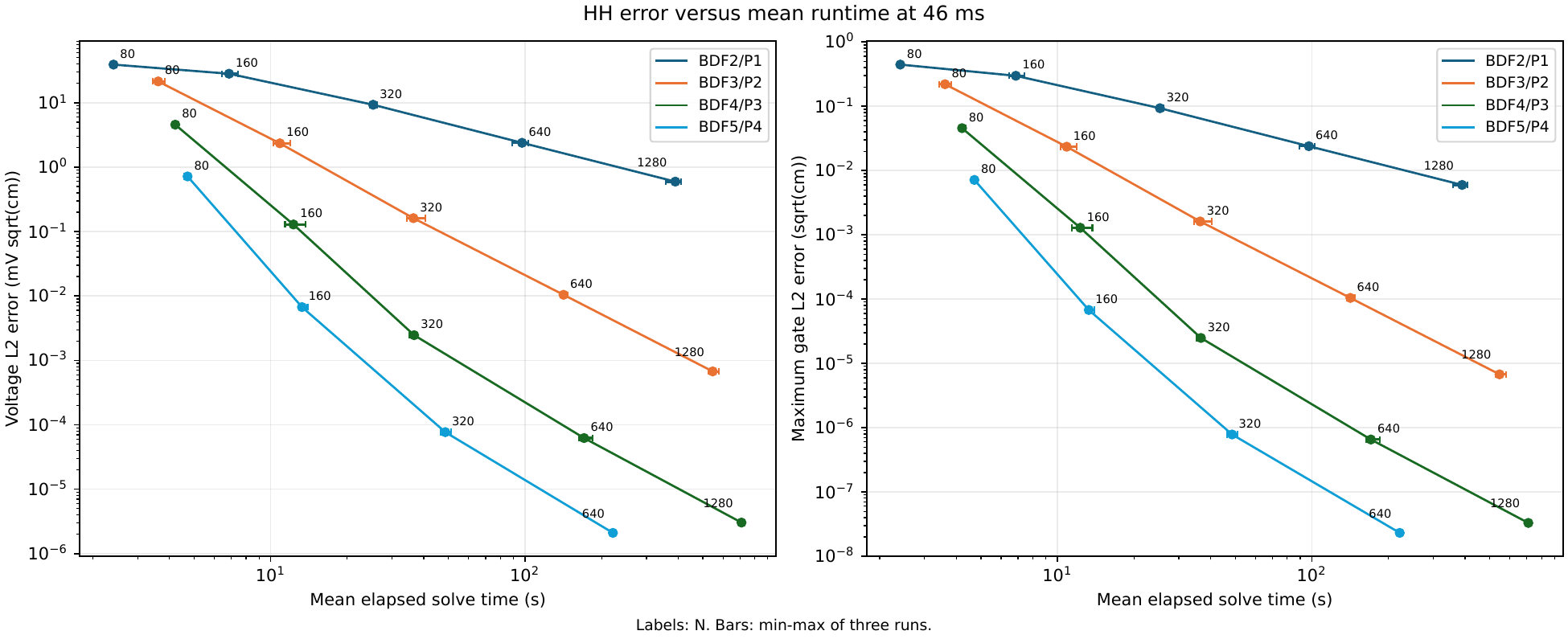}
		\caption{Final-time $L^2(\Gamma)$ errors versus mean elapsed runtime
		at $T=46\,\mathrm{ms}$: voltage (left) and maximum gate error (right).
		Point labels give $N_e$; horizontal bars show the range of three runs.}
		\label{fig:runtime_benchmark}
	\end{figure}

	Each of the $19$ cases is computed three times in fresh serial processes,
	with an identical untimed warm-up and randomized case order within each
	repetition block.  All $57$ runs use one numerical-library thread on an
	Intel Core i5-12400, with Python~3.13.13, NumPy~2.4.4, and SciPy~1.17.1.
	Elapsed wall-clock time includes mesh construction, solver setup,
	starting histories, and time integration.  Imports, warm-up, reference
	loading, file output, and error evaluation are excluded.
	Figure~\ref{fig:runtime_benchmark} plots the arithmetic mean of the three
	runtimes against the final-time errors; horizontal bars show the
	minimum and maximum runtimes.

	For voltage $L^2(\Gamma)$ error $E_V^{46}\le10^{-3}$, the fastest tested
	cases are listed in Table~\ref{tab:runtime_accuracy_target}.
	BDF2--$P_1$ does not reach this target on the tested grids.
	BDF5--$P_4$ has the lowest mean runtime at this target among the tested
	joint refinement paths.  All errors are relative to the same finite
	numerical reference.

	\begin{table}[htbp]
		\centering
		\caption{Fastest tested cases with $E_V^{46}\le10^{-3}$.}
		\label{tab:runtime_accuracy_target}
		\begin{tabular}{@{}l r r r@{}}
			\toprule
			Scheme & $N_e$ & $E_V^{46}$ & Mean elapsed time (s) \\
			\midrule
			BDF3--$P_2$ & $1280$ & $6.7427\times10^{-4}$ & $546.7$ \\
			BDF4--$P_3$ & $640$ & $6.2218\times10^{-5}$ & $170.6$ \\
			BDF5--$P_4$ & $320$ & $7.6875\times10^{-5}$ & $48.6$ \\
			\bottomrule
		\end{tabular}
	\end{table}

	\section{Conclusions}
	
	We have developed and analyzed a linearly implicit IMEX BDF--FEM family
	for the Hodgkin--Huxley cable system on a branched metric graph, using
	constant time steps, orders $q=1,\ldots,5$, and any fixed spatial degree
	$p\ge1$.  Smooth masking of the extrapolated gates preserves
	voltage-equation coercivity without projecting the numerical gates, and
	a voltage cutoff gives bounded Lipschitz gate rates.  Each timestep
	requires one linear voltage solve followed by independent nodal gate
	updates.  Under exact spatial integration, the complete regularized scheme
	satisfies unconditional bounds in discrete history energies.  Under
	regularity and starting-history assumptions, the voltage and reconstructed
	gates converge in $L^2(\Gamma)$ at the optimal rate
	$O(h^{p+1}+\Delta t^q)$, without a timestep restriction or an
	$h$--$\Delta t$ coupling.
	
	Manufactured solutions recover the predicted orders through BDF5--$P_4$
	on the Y-shaped graph, and an additional two-junction test verifies
	assembly and spatial accuracy with unequal edge lengths.  Physical
	simulations demonstrate pulse propagation and splitting on the Y-shaped
	axon.  In the controlled work--precision comparison within the proposed
	family, higher-order variants attain the tested accuracy target at lower
	runtime, with BDF5--$P_4$ fastest among the tested cases.  These comparisons
	use a common finite numerical reference and the tested joint-refinement
	paths.
	
	For numerical quadrature, the convergence rate remains conditional on
	the stated bound on the full exact-form residual.  The Gauss--Legendre
	implementation exhibits the expected rates, while a rule-specific proof
	of this defect bound for the masked coefficients and stimulus remains
	open.  Such an analysis should identify practical quadrature orders that
	preserve the theoretical accuracy uniformly under refinement.
	
	Under the stated assumptions, the exact-integration analysis extends to
	any fixed connected finite metric graph with standard continuity and
	Kirchhoff conditions, sealed outer endpoints, and uniform cable and gate
	laws.  The experiments use uniform coefficients and smooth data; physical
	pulse propagation is tested on a symmetric Y-shaped graph, while the
	two-junction manufactured fields have zero edgewise derivatives at every
	vertex.  Testing pulses through successive junctions and extending the
	method to heterogeneous networks are natural next steps.  Spatial and
	temporal adaptivity could improve efficiency when regularity varies, but
	variable time stepping requires further stability and convergence
	analysis.  Solver and preconditioning studies should also address
	scalability on larger neuronal graphs.

	\bibliographystyle{plain}
	\bibliography{ref}

	\clearpage
	\appendix
	\section{The $\mathcal{C}^{2N}$ Continuous Masking Function}\label{app:masking_function}
	For an integer $N\ge1$, define the truncated Maclaurin series
	of $\operatorname{arctanh}(u)$ by
	\begin{equation}
		P_N(u) = \sum_{k=1}^N \frac{u^{2k-1}}{2k-1}.
	\end{equation}
	
	For $\varepsilon>0$, set
	\begin{equation}
		\phi(x) =
		\begin{cases}
			\varepsilon \tanh\left( P_N\left(\dfrac{x}{\varepsilon}\right) \right), & \text{if } x < 0, \\[10pt]
			x, & \text{if } 0 \le x \le 1, \\[10pt]
			1 + \varepsilon \tanh\left( P_N\left(\dfrac{x-1}{\varepsilon}\right) \right), & \text{if } x > 1.
		\end{cases}
	\end{equation}
	This function equals the identity on $[0,1]$ and takes values in
	$[-\varepsilon,1+\varepsilon]$.

	\subsection{Smoothness}
	Since $P_N(u)$ agrees with $\operatorname{arctanh}(u)$ through degree
	$2N-1$, we have
	\[
	\tanh(P_N(u))=u+\mathcal{O}(u^{2N+1})
	\qquad \text{as } u\to0.
	\]
	Thus the Taylor series around each patching point matches the identity
	function through order $2N$, and the piecewise function is $\mathcal{C}^{2N}$.
	
	\subsection{Global Lipschitz Constant}
	On either outer branch, the derivative has the form
	\[
	\phi'(x)=\frac{P_N'(u)}{\cosh^2(P_N(u))},
	\qquad P_N'(u)=\sum_{k=0}^{N-1}u^{2k}>0,
	\]
	where $u=x/\varepsilon$ or $u=(x-1)/\varepsilon$.
	Since $P_N(u)-\operatorname{arctanh}(u)=\mathcal O(u^{2N+1})$ and
	$(\cosh^2)'(z)=\sinh(2z)=\mathcal O(z)$ near zero,
	\[
	\cosh^2(P_N(u))=\frac1{1-u^2}+\mathcal O(u^{2N+2})
	\qquad\text{as }u\to0.
	\]
	The entire function $\cosh^2(P_N(u))$ has a globally convergent even
	Taylor series with nonnegative coefficients: $P_N$ has nonnegative
	odd coefficients, and $\cosh^2 z$ has nonnegative even coefficients.
	The local expansion identifies the coefficients through degree $2N-2$
	as one.  Therefore, for every real $u$,
	\[
	\cosh^2(P_N(u))
	=\sum_{k=0}^{N-1}u^{2k}+\sum_{k=N}^{\infty}c_ku^{2k}
	\ge P_N'(u),
	\qquad c_k\ge0.
	\]
	Thus $0<\phi'(x)\le1$ globally.  Since $\phi'(x)=1$ on $(0,1)$,
	the optimal Lipschitz constant is $L_\phi=1$.

	\subsection{The $C^6$ Mask Used in the Computations}
\begin{figure}[htp]
		\centering
		\includegraphics[width=0.6\textwidth]{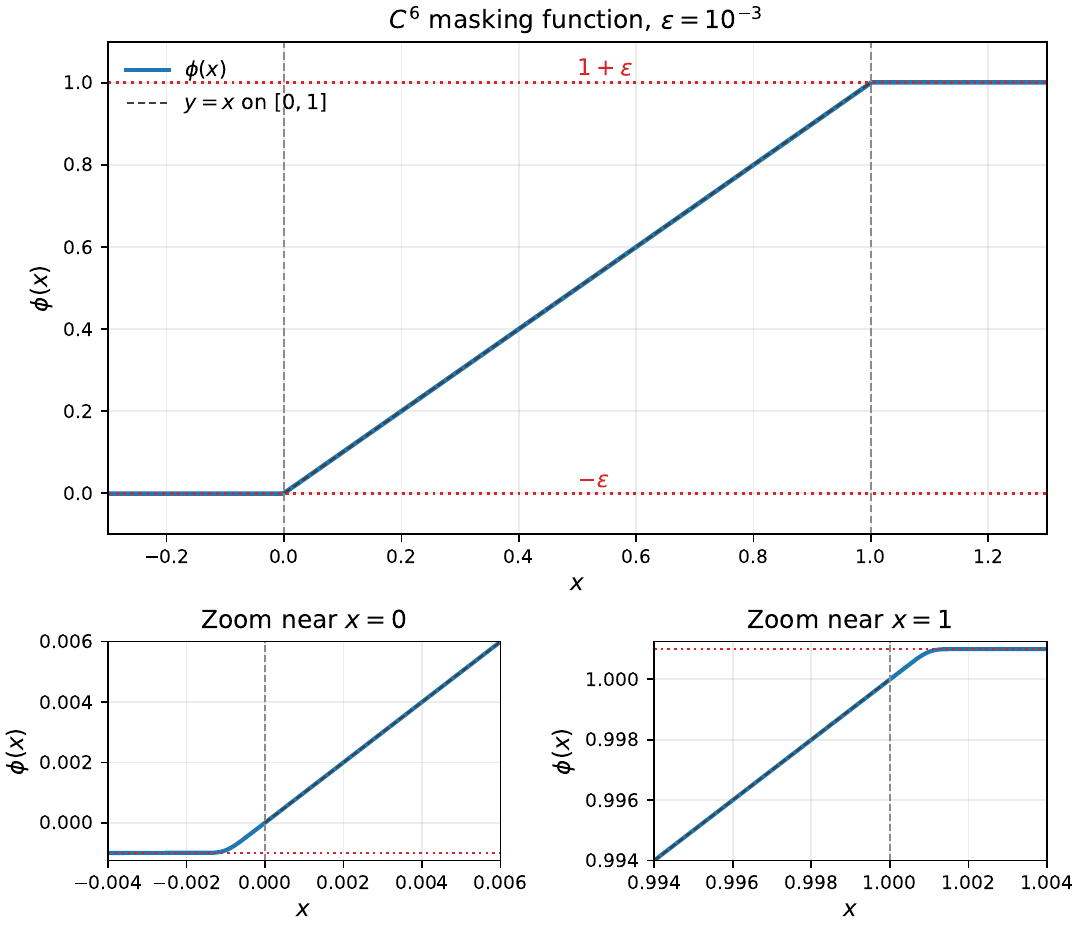}
		\caption{The $\mathcal C^6$ masking function with $N=3$ and
		$\varepsilon=10^{-3}$, including details near the patching points.}
		\label{fig:masking_function}
	\end{figure}

	The choice $N=3$ gives
	\begin{equation}
		P_3(u)=u+\frac{u^3}{3}+\frac{u^5}{5}.
	\end{equation}
	The preceding arguments establish $\phi\in\mathcal C^6(\mathbb R)$
	and $L_\phi=1$.  More precisely, as $x\to0^-$,
	\[
	\varepsilon\tanh\left(P_3\left(\frac{x}{\varepsilon}\right)\right)
	=x-\frac{x^7}{7\varepsilon^6}
	+\mathcal O\left(\frac{x^9}{\varepsilon^8}\right).
	\]
	The nonzero seventh-order term prevents $\mathcal C^7$ matching with
	the identity at $x=0$; the same conclusion holds at $x=1$ by translation.
	Figure~\ref{fig:masking_function} shows this mask for the value of
	$\varepsilon$ used in the experiments.

\end{document}